\documentclass[12pt]{amsart}
\usepackage[dvips]{graphicx,color}
\usepackage{amssymb,amscd,latexsym,epsfig,enumerate}
\usepackage[mathscr]{euscript}
\usepackage[curve]{xypic}
\usepackage{setspace}

\DeclareRobustCommand{\SkipTocEntry}[4]{}
\makeatletter
\newcommand\@dotsep{4.5}
\def\@tocline#1#2#3#4#5#6#7{\relax
  \ifnum #1>\c@tocdepth % then omit
  \else
    \par \addpenalty\@secpenalty\addvspace{#2}%
    \begingroup \hyphenpenalty\@M
    \@ifempty{#4}{%
      \@tempdima\csname r@tocindent\number#1\endcsname\relax
    }{%
      \@tempdima#4\relax
    }%
    \parindent\z@ \leftskip#3\relax \advance\leftskip\@tempdima\relax
    \rightskip\@pnumwidth plus1em \parfillskip-\@pnumwidth
    #5\leavevmode\hskip-\@tempdima #6\relax
    \leaders\hbox{$\m@th
      \mkern \@dotsep mu\hbox{.}\mkern \@dotsep mu$}\hfill
    \hbox to\@pnumwidth{\@tocpagenum{#7}}\par
    \nobreak
    \endgroup
  \fi}
\makeatother 

\DeclareFontFamily{OT1}{rsfs}{}
\DeclareFontShape{OT1}{rsfs}{n}{it}{<-> rsfs10}{}
\DeclareMathAlphabet{\curly}{OT1}{rsfs}{n}{it}

\newcommand\I{\curly I}

\newcommand\LL{\mathbb L}
\renewcommand\O{\mathcal O}
\newcommand\oh{\mathcal O}
\newcommand\PP{\mathbb P}

\newcommand\M{\mathcal M}
\newcommand\MM{\,\overline{\!\mathcal{M}}}
\newcommand\C{\mathbb C}

\newcommand\Z{\mathbb Z}
\newcommand\com{{\mathbb{C}}}
\newcommand\rarr{\rightarrow}
\newcommand\proj{{\mathbb{P}}}

\newcommand{\rt}[1]{\stackrel{#1\,}{\rightarrow}}
\newcommand{\Rt}[1]{\stackrel{#1\,}{\longrightarrow}}
\newcommand\To{\longrightarrow}
\newcommand\into{\hookrightarrow}
\newcommand\Into{\ar@{^{ (}->}[r]}

\newcommand\bull{{\scriptscriptstyle\bullet}}
\newcommand\udot{^\bull}
\newcommand\dbar{\overline\partial}

\newcommand\QQ{{\mathcal{Q}}}
\newcommand\rk{\operatorname{rank}}

\newcommand\Hom{\operatorname{Hom}}

\newcommand\Ext{\operatorname{Ext}}
\newcommand\ext{\curly Ext}

\newcommand\Pic{\operatorname{Pic}}

\newcommand\Sym{\operatorname{Sym}}

\newcommand\beq[1]{\begin{equation}\label{#1}}
\newcommand\eeq{\end{equation}}
\newcommand\beqa{\begin{eqnarray*}}
\newcommand\eeqa{\end{eqnarray*}}

\newcommand\halfsection[1]{\renewcommand{\thesection}{$\frac12$}
\section{#1}\renewcommand{\thesection}{\arabic{section}}}

\newcommand\newsection[1]{\renewcommand{\thesection}{\arabic{section}$\frac12$}
\section{#1}\renewcommand{\thesection}{\arabic{section}}}

\makeatletter \@addtoreset{equation}{section} \makeatother

\title{13/2 ways of counting curves}
\author[R. Pandharipande, R. P. Thomas]{R. Pandharipande and R. P. Thomas}

\begin{document}

{\tiny
\begin{tabbing}
{\ \ \ \ \ \ \ \ \ \ \ \ \ \ \  \ \ \ \  
\ \ \ \ \ \ \ \ \  \ \ \ \  \ \ 
\ \ \ \ \ \ \ \ \  \ \ \ \  \ \ 
 \ \ \ \ \ \ \ \ \ \ \ \ \ \  \ \ }
\={\em I was of three minds},\\  \>{\em Like a tree} \\
\>{\em In which there are three blackbirds.} \\ \\
\> \ \ \ \ \ \ \ \ \ \ \ \ \ Wallace Stevens
\end{tabbing}
}

\vspace{30pt}

\begin{abstract} \noindent
In the past 20 years, compactifications of the 
families of curves in algebraic
varieties $X$ have been studied via stable maps, Hilbert schemes, stable
pairs, unramified maps, and stable quotients. Each path leads to   
a different enumeration of curves.
A common thread is the use of a 2-term deformation/obstruction theory
to define a virtual fundamental class.
The richest geometry occurs when  $X$ is a nonsingular projective variety
of dimension 3.

We survey here the 13/2 principal ways to count curves with special
attention to the 3-fold case. The different theories are 
linked by a web of conjectural relationships which we highlight. 
Our goal is to provide a guide for 
graduate students looking for an elementary route 
into the subject.
\end{abstract}

\maketitle

\begin{spacing}{1.3}
\tableofcontents
\end{spacing}

\pagebreak

%%%%%%%%%%%%%%%%%%%%%%%%%%%%%%%%%%%%%%%%%%%%%%%%%%%%%%%%%%%%%%%%%%%%%%

\setcounter{section}{-1}
\section{Introduction}

\subsection*{\S Counting}
Let $X$ 
be a nonsingular projective variety (over $\mathbb{C}$), and let
$\beta\in H_2(X,\Z)$ be a homology class.
We are interested here in
 counting the  algebraic curves of $X$ in class $\beta$.
For example, how many twisted cubics in $\PP^3$ meet
 12 given lines? Mathematicians such as Hurwitz, Schubert, and Zeuthen 
 have considered such questions since the $19^{th}$
century. Towards the end of the $20^{th}$ century and continuing to
the present,  the subject has been greatly enriched by
 new insights from symplectic geometry and 
topological string theory. 

Under appropriate genericity conditions,
counting solutions in algebraic geometry often yields 
deformation invariant answers. 
A simple example is provided by Bezout's Theorem
concerning the intersections of plane curves. 
Two generic algebraic
curves in $\C^2$
of degrees $d_1$ and $d_2$ intersect transversally in finitely many points. Counting these points
yields the topological intersection number $d_1d_2$. But in nongeneric situations,
 we can find fewer solutions or an infinite number. 
The curves may intersect with tangencies in a smaller number of points (remedied by counting intersection points with multiplicities).
If the curves intersect ``at infinity", we will again find fewer intersection points in $\C^2$
whose total we do not consider to be a ``sensible'' answer. 
Instead, we compactify $\C^2$ by $\PP^2$ and count there. Finally, the curves may intersect in an entire component.
The technique of excess intersection theory is  required then to  
obtain the correct answer. 
Compactification and transversality already play a important
role in the geometry of Bezout's Theorem.

Having deformation invariant answers for the enumerative
geometry of curves in $X$ is desirable for several reasons. 
The most basic is the possibility of deforming $X$ to
a more convenient space.
To achieve deformation invariance, two main issues must be considered:
\begin{enumerate}
\item[(i)] compactification of the moduli space $\M(X,\beta)$ of curves \\
$C\subset X$ of class $\beta$,
\item[(ii)] transversality of the solutions.
\end{enumerate}
What we mean by the moduli space $\M(X,\beta)$ is to be explained
 and will differ in each of the sections below.
Transversality concerns both the possible excess dimension
of $\M(X,\beta)$ and the transversality of the constraints.

\subsection*{\S Compactness}

For Bezout's Theorem, 
we compactify the geometry so 
intersection points running to infinity do not escape 
our counting. The result is a deformation invariant answer. 

A compact space $\M(X,\beta)$ 
which parameterises all nonsingular embedded curves in class $\beta$ 
will usually have to contain singular curves of some 
sort.  Strictly speaking, the compact moduli spaces $\M(X,\beta)$
will often not be compactifications of  
the spaces of nonsingular embedded curves
--- the latter need not be dense in $\M(X,\beta)$. For instance 
$\M(X,\beta)$ might be nonempty when there are 
{\em no} nonsingular embedded curves.
The singular strata are important for deformation 
invariance. As we deform $X$, curves can ``wander off to infinity" in $\M(X,\beta)$ by becoming singular.

\subsection*{\S Transversality}
A simple question to consider is the number of 
elliptic cubics in $\PP^2$ passing through 9 points $p_1,\ldots,p_9\in
\PP^2$.
The linear system 
$$\PP(H^0(\PP^2,\O_{\PP^2}(3)))\cong\PP^9$$
provides a natural compactification of the moduli space.
%To cut down to a 0-dimensional space  whose points we 
%can try to count, fix 9 points $p_i$ in $\PP^2$ and consider 
%those curves that pass through them.
%Equivalently, 
%blow $\PP^2$ up in the $p_i$ and consider curves in the 
%class of the pullback of $\O(3)$ minus the 9 exceptional divisors.
Each $p_i$ imposes a single linear condition which
determines a hyperplane $$\PP^8_i\subset\PP^9,$$ 
of curves passing through $p_i$. For general $p_i$, these 9 
hyperplanes are transverse and intersect in a single point. 
Hence, we expect our count to be 1. 
But if the $p_i$ are the 9 intersection points of two cubics,
then we obtain an entire pencil of solutions given by
the linear combinations of the two cubics.

An alternative way of looking at the same enumerative question is the
following.
Let $$\epsilon:S \rightarrow \PP^2$$
be the   
blow-up of $\PP^2$  at 9 points  $p_i$ and consider curves in the 
class 
$$\beta= 3H -E_1-E_2- \ldots - E_9$$
where $H$ is the $\epsilon$ pull-back of the hyperplane class and
 the $E_i$ are the exceptional divisors. In general there will be a unique elliptic curve embedded in class $\beta$. But
if the 9 points are the intersection of two cubics, then $S$
is a rational elliptic surface via the pencil
$$\pi: S \rightarrow \PP^1 .$$ 
How to sensibly ``count" the pencil of elliptic fibres on $S$
is not obvious.

A temptation based
on the above discussion is to define the enumeration of curves by
counting {after taking a generic perturbation of the geometry}. 
Unfortunately, we often
do not have enough perturbations to make the situation
fully transverse. A basic rigid example 
 is given by counting the intersection points of a $(-1)$-curve with itself
on a surface.
Though we cannot algebraically
move the curve to be transverse to itself,
we know another way to get the ``sensible" answer of topology:
take the Euler number $-1$ of the normal bundle. 
In curve counting, there is a similar excess intersection theory 
approach to getting a sensible, deformation invariant answer using 
{\em virtual fundamental classes}.

For the rational elliptic surface $S$,
 the base $\PP^1$ is
a natural compact moduli space 
parameterising the elliptic curves in the pencil. 
%The reader tempted to try to define invariants using only the smooth curves in %the pencil should remember that although generically there are 12 singular (nod%al) curves in the pencil, one can deform to smooth surfaces with fewer singular% fibres (with worse singularities).
The count of elliptic fibres is the Euler class of the {\em
obstruction bundle} over the pencil $\PP^1$. Calculating the obstruction bundle to be $\O_{\PP^1}(1)$,
we recover the answer 1 expected from deformation invariance.

Why is the obstruction bundle $\O_{\PP^1}(1)$?
In Section \ref{GW}, a short introduction to the deformation
theory of maps is presented. Let
$E\subset S$
be the fibre of $\pi$ over $[E]\in \PP^1$.
Let $\nu_E$ be the normal bundle of $E$ in $S$.
The obstruction space at
$[E]\in \PP^1$
is 
$$H^1(E,\nu_E) = H^1(E,\O_E) \otimes \O_{\PP^1}(2)|_{[E]}\ .$$
The term $H^1(E,\O_E)$ yields the dual of the Hodge bundle
as $E$ varies 
and is isomorphic to $\O_{\PP^1}(-1)$.
Hence, we find the obstruction bundle to be
$\O_{\PP^1}(1)$.

We will discuss virtual classes in the Appendix.  
We should think loosely of $\M(X,\beta)$ as being cut out of 
a nonsingular ambient space by a set of equations. 
The expected, or {\em virtual}, dimension of $\M(X,\beta)$ is the dimension of the ambient space minus the number of equations. If the
derivatives of the equations are not linearly independent 
along the common zero locus $\M(X,\beta)$,
 then $\M(X,\beta)$ will be singular or have dimension higher than expected.
In practice, $\M(X,\beta)$ is very rarely nonsingular of the
expected dimension. 
We should think of the virtual class as representing the fundamental
cycle of the  ``correct" moduli space 
(of dimension \emph{equal} to the virtual dimension) inside the actual moduli space. The virtual class may be considered to 
give the result of perturbing the setup to a transverse geometry, 
even when such perturbations do not actually exist.

\subsection*{\S Overview}
A nonsingular embedded curve $C\subset X$ can be described in two 
fundamentally different ways:
\begin{enumerate}
\item[(i)] as an algebraic map $C\to X$
\item[(ii)] as the zero locus of an ideal of algebraic functions on $X$.
\end{enumerate}
In other words, $C$ can be seen as
a \emph{parameterised} curve with a map or
an \emph{unparameterised} curve with an embedding.
Both realisations
arise naturally in physics --- the first as the worldsheet of a string moving
in $X$, the second as a D-brane or boundary condition embedded in $X$.

Associated to the two basic ways of thinking of curves, there are two natural
paths for compactifications. The first 
allows the map $f$ to degenerate badly while keeping the domain
curve as nice as possible.
%
%leading to stable maps and Gromov-Witten theory, as described in Section \ref{GW}. 
The second
keeps the map as an embedding but allows the curve to degenerate arbitrarily.

%
%(to any 1-dimensional subscheme of $X$) and leads to the Hilbert scheme compactification, of Section \ref{MNOPsec}.

We describe here $6\frac{1}{2}$ 
methods for defining curve counts in
algebraic geometry. We start in Section $\frac{1}{2}$ with a
discussion of the successes
and limitations of the naive counts pursued by the $19^{th}$ century
geometers (and followed for more than 100 years). Since 
such counting is not always well-defined and
has many drawbacks, we view the naive approach as only $\frac{1}{2}$
a method.

In Sections $1\frac{1}{2}$ -- $6\frac{1}{2}$,
six approaches to deformation invariant curve counting
are presented.
Two (stable maps and unramified maps) fall in class (i),
three (BPS invariants, ideal sheaves, stable pairs) in class (ii),
and one (stable quotients) straddles both classes (i-ii). 
The compactifications and 
virtual class constructions are dealt with 
differently in the six cases. Of course,
each  of the six has advantages and drawbacks.

There are several excellent references covering 
different aspects of the material surveyed here in much
greater depth, see for instance \cite{CKK, Hori,
LiSurvey, PandClay, TodaSurvey}. 
Also, there are many
beautiful directions which
we do not cover at all. For example,  mirror symmetry, integrable hierarchies,
 descendent
invariants, 3-dimensional
partitions,
 and holomorphic symplectic geometry 
all play significant roles in the subject. Though
orbifold and relative geometries
have been very important for the development of the ideas
presented here, we have chosen to omit a discussion.
Our goal is to describe the $6\frac12$ counting theories
as simply as possible and
to present the web of relationships amongst them.

%The resulting two sets of invariants have their advantages and drawbacks, which we outline. The other sections describe variants that overcome some of these drawbacks. The stable pairs of Section \ref{PT} use a compactification closely related to the Hilbert scheme of Section \ref{MNOPsec} which is in some natural sense smaller and easier to relate to curve counting. The unramified maps of Kim et al (Section \ref{Kim}) do a similar thing for the stable maps of Gromov-Witten theory. The stable quotients of Section \ref{MOP} share features of both types of compactifications. Finally there are the BPS invariants of Gopakumar-Vafa. These first arose as an attempt to describe Gromov-Witten theory via the physics of M-theory, involving a sheaf theory more closely related to the second type of compactification. In many ways they should be the ideal curve counting invariants (at least for 3-folds), promising to combine the advantages of each of the other theories. For now they are still conjectural in mathematics, but their link to stable pairs has led to a lot of recent progress which we survey briefly in Section \ref{GVsec}.

\subsection*{\S Acknowledgments}
We thank our students, collaborators, and
colleagues for their contributions to our
understanding of the subjects presented here.
Special thanks are due to L. Brambila-Paz,
J. Bryan, Y. Cooper, S. Katz, A. Kresch, A. MacPherson, P. Newstead and H.-H.
Tseng for their 
specific comments and suggestions
about the paper. 

R.P. was partially supported by NSF grant 
DMS-1001154, a Marie Curie fellowship
at 
IST Lisbon, and 
a grant from the Gulbenkian foundation.
R.T. was partially supported by an EPSRC programme grant. We would both like to thank the Isaac Newton Institute, Cambridge for support and a great research environment.

\halfsection{Naive counting of curves} 
Let $X$ be a nonsingular projective variety, and let 
$\beta\in H_2(X,\mathbb{Z})$ be a homology class.
Let 
$C\subset X$
be a nonsingular embedded (or immersed) curve of genus $g$ and class
$\beta$. 
The expected dimension of the family of genus $g$ and class $\beta$
curves containing $C$ is
\begin{equation}\label{exx}
3g-3+ \chi(T_X|_C) = \int_C c_1(X)+ ( {\text{dim}}_\com X -3)(1-g)\,.
\end{equation}
The first  term on the left comes from the complex moduli of the 
genus $g$ curve,
$$\text{dim}_\com \ \MM_g = 3g-3  \,.$$
The second term  arises from
infinitesimal
deformations of $C$ which do not change the complex structure of $C$.
More precisely,
$$\chi(T_X|_C)= h^0(C, T_X|_C)- h^1(C,T_X|_C)$$
where  $H^0(C, T_X|_C)$ is the space of such deformations (at least
when $C$ has no continuous families of automorphisms).
The ``expectation'' amounts to the vanishing of
 $H^1(C, T_X|_C)$. Indeed if $H^1(C, T_X|_C)$ vanishes, the
family of curves is nonsingular of expected dimension at $C$, see
 \cite{Kollar}.
%More generally, as discussed in the Appendix, the implicit function theorem and Kuranishi theory imply that the family of curves with fixed complex structure in $X$ can be written locally as a subscheme of the tangent space $H^0(C, T_X|_C)$, cut out by an equation with values in $H^1(C, T_X|_C)$. 
We will return to this deformation theory in Section \ref{GW}.

If the open family of embedded (or immersed) curves of genus $g$ 
and class $\beta$
is of pure expected dimension \eqref{exx}, then naive classical
curve counting is sensible to undertake. We can attempt to count the
actual  numbers of embedded (or immersed) curves
of genus $g$ and class $\beta$ in $X$ subject
to incidence conditions.

\bigskip

The main classical{\footnote{We
do not attempt here to give a complete
classical bibliography. Rather, the references
we list, for the most part, are modern
treatments.}} examples where naive curve counting with
simple incidence is 
reasonable to consider constitute a rather short list:
\vspace{10pt}
\begin{enumerate}
\item[(i)] Counting Hurwitz coverings of $\PP^1$ and curves of higher genus,
\item[(ii)] Severi degrees in $\PP^2$ and $\PP^1 \times \PP^1$ in all genera,
\item[(iii)] Counting genus 0 curves in general blow-ups of $\PP^2$,
\item[(iv)] Counting genus 0 curves in homogeneous spaces
 such as $\PP^n$, Grassmannians, and flag varieties,
\item[(v)] Counting lines on complete intersections in $\PP^n$,
\item[(vi)] Counting curves of genus 1 and 2 in $\PP^3$.
\end{enumerate}
\vspace{10pt}
The Hurwitz covers of $\PP^1$ (or higher genus curves), 
$$C\rightarrow \PP^1,$$ 
are neither embeddings
nor immersions, but rather are counts of ramified maps, see
\cite{OP1} for an introduction. Nevertheless (i) fits naturally
in the list of classical examples. The Severi degrees (ii) are
the numbers of immersed curves of genus $g$ and class $\beta$
passing through the expected number of points 
on a surface. Particularly for the case of $\PP^2$, the study 
of Severi degrees has a long
history \cite{CH,Gott,Ziv}. Counting genus 0 curves  on
blow-ups (iii) is equivalent to imposing multiple point singularities
for plane curves, 
see \cite{GotP} for a treatment.
Genus 0 curves behave very well in homogeneous
spaces, so the questions (iv) have been considered since
Schubert and Zeuthen \cite{Sch, Zeu}.
Examples of (v) include the famous 27 lines on a cubic
surface and the 2875 lines on a quintic 3-fold, see \cite{CKK}. The 
genus 1 and 2 enumerative geometry of space curves 
was much less studied by the classical geometers,
but still can be viewed in terms of  naive counting.

For particular genera and classes on other varieties
$X$, the families of curves might be pure of expected
dimension. The above list addresses the cases of more uniform
behavior.
Until new ideas from symplectic
geometry and topological string theory
were introduced in the 1980s and 90s, the
classical cases (i-vi)
were the main
topics of study in enumerative geometry.
The subject was an important
area, especially for the development of
intersection theory in algebraic geometry.
See \cite{Fu} for a historical survey.
However,
because of the restrictions,
we treat naive counting as only $\frac{1}{2}$ of
an enumerative theory here.

New approaches to
enumerative geometry by tropical methods
have been developed extensively in recent years \cite{IMS,Mik}.
However, the lack of a virtual fundamental
class in tropical geometry restricts the 
direct{\footnote{Tropical methods
do interact in an intricate way with
virtual curve counts on Calabi-Yau 3-folds in
the program of Gross and Siebert \cite{Grr,GPS, GrS2} to study mirror symmetry.}}
 applications at the moment to the classical cases.

The counting of rational curves on algebraic $K3$ surfaces is almost
a classical question. A $K3$
surface with Picard number 1 has finitely many rational
curves in the primitive class (even though the expected
dimension of the family of rational curves is $-1$
 by \eqref{exx}). As proved in \cite{XC}, for a general
$K3$ of Picard number 1, all the primitive rational curves
are nodal.
A proposal for the
count was made by Yau and Zaslow \cite{YZ} in terms of modular
forms. The proofs by Beauville \cite{Beau} and Bryan-Leung \cite{BryL} certainly
use modern methods. The counting of rational curves in all
(including imprimitive) classes on $K3$ surfaces
shows the fully non-classical
nature of the question \cite{NLYZ}.

\setcounter{section}{0}

\newsection{Gromov-Witten theory} \label{GW}

\subsection*{\S Moduli} Gromov-Witten theory
provided the first modern approach
to curve counting which dealt successfully with the
issues of compactification and transversality. 
The subject has origins in Gromov's work on pseudo-holomorphic
curves in symplectic geometry \cite{Gr} and papers of Witten on 
topological strings \cite{Wit}.
Contributions by Kontsevich, Manin, Ruan, and Tian
 \cite{K2,KoMa, Ruan, RuTian} 
played an important
role in the early development.

In Gromov-Witten theory, curves are viewed as parameterised with 
an algebraic map
$$C \rightarrow X\,. $$
The compactification strategy is to 
admit only nodal singularities in the domain 
while allowing the map to become rather degenerate. 
More precisely, define $\MM_g(X,\beta)$ to be the moduli space of 
{\em stable maps}:
$$
\left\{f\colon C\to X\,\Bigg|
\begin{array}{c}C\mathrm{\ a\ nodal\ curve \ of\ arithmetic\ genus\ }g,\\ 
f_*[C]=\beta, \mathrm{\ and\ Aut}(f)
\mathrm{\ finite} \end{array}\!\!\!\right\} .
$$
The map $f$ is invariant under an
automorphism $\phi$ of the domain $C$ if
$$f= f\circ \phi\,.$$
By definition, $\mathrm{Aut}(f)\subset \mathrm{Aut}(C)$ 
is the subgroup of elements for which $f$ is invariant.
The finite automorphism condition for a stable map implies 
the moduli space $\MM_g(X,\beta)$ is naturally a
Deligne-Mumford stack.

The compactness of $\MM_g(X,\beta)$ is not immediate. A proof
can be found in \cite{FuPa} using standard properties of
semistable reduction for curves. In Section \ref{MNOPsec} below,
we will discuss nontrivial limits in the space of stable maps,
see for instance \eqref{conics} and \eqref{genuschange}.

\subsection*{\S Deformation theory} We return now to the deformation
theory for embedded curves 
briefly discussed  in 
Section $\frac{1}{2}$. The deformation theory for arbitrary stable
maps is very similar.

Let $C\subset X$ be a nonsingular embedded curve
with normal bundle $\nu\_C$.
 The Zariski tangent space to the moduli space $\MM_g(X,\beta)$ 
at the point $[C\rightarrow X]$ is 
given by $H^0(C,\nu\_C)$. Locally, we can lift 
a section of $\nu\_C$ to a section of $T_X|_C$ and deform $C$ along 
the lift to first order. Since globally $\nu\_C$ is not usually a summand of 
$T_X|_C$ but only a quotient, the lifts will differ over overlaps 
by vector fields along $C$. The deformed curve will have a 
complex structure whose transition functions differ by these vector fields. 
In other words, from
$$
0\to T_C\to T_X|_C\to\nu\_C\to0 ,
$$
we obtain the sequence
\begin{equation}\label{hj34}
0\to H^0(C,T_C)\to H^0(C,T_X|_C)\to H^0(C,\nu\_C)\to H^1(C,T_C)
\end{equation}
which expresses how deformations in $H^0(C,\nu\_C)$ change 
the complex structure on $C$ through the boundary map to $H^1(C,T_C)$. 
The kernel $$H^0(C,T_X|_C)/H^0(C,T_C)$$ consists of the deformations 
given by moving $C$ along vector fields in $X$, thus preserving 
the complex structure of $C$, modulo infinitesimal automorphisms of $C$. 
Similarly, obstructions to deformations lie in $H^1(C,\nu\_C)$.

The expected dimension $\chi(\nu\_C)=h^0(\nu\_C)-h^1(\nu\_C)$ of the moduli space is given by the calculation
\begin{equation}\label{hj67}
\chi(\nu\_C) = \int_C c_1(X) + (\text{dim}_\com X-3)(1-g)
\end{equation}
obtained from sequence \eqref{hj34}. If $H^1(C, T_X|_C)$
vanishes, so does the obstruction space
$H^1(C,\nu\_C)$. Formula \eqref{hj67} then
computes the actual dimension of the Zariski tangent space.

For arbitrary stable maps $f:C \rightarrow X$,
we replace the dual of $\nu\_C$ by the complex 
\begin{equation}
\{f^*\Omega_X\to\Omega_C\} \label{rt45}
\end{equation}
 on $C$. 
If $C$ is nonsingular and $f$ is an embedding, the complex \eqref{rt45}
is quasi-isomorphic to its kernel $\nu_C^*$. 
The deformations/obstructions of $f$ are governed by
\beq{GWdef}
\Ext^i\Big(\{f^*\Omega_X\to\Omega_C\},\O_C\Big)
\eeq
for $i=0,1$. Similarly the deformations/obstructions of $f$ with the curve $C$ fixed are governed by $\Ext^i\big(f^*\Omega_X,\O_C\big)=H^i(f^*T_X)$.

Since the Ext groups \eqref{GWdef} vanish for $i\ne0,1$, the 
deformation/ob\-struction theory is 2-term.
The moduli space 
admits a virtual fundamental class{\footnote{The virtual fundamental class
is algebraic, so should be more naturally considered in the Chow group
$A_*(\MM_g(X,\beta),\mathbb{Q})$.}} 
$$[\MM_g(X,\beta)]^{vir} \in H_*(\MM_g(X,\beta),\mathbb{Q})$$
of complex dimension equal to the virtual dimension
\beq{vd}
\text{ext}^0-\text{ext}^1=\int_\beta c_1(X)+(\dim_\com X-3)(1-g)\,.
\eeq
An introduction to the virtual fundamental
class is provided in the Appendix.

\subsection*{\S Invariants}
To obtain numerical invariants, we must cut the virtual class from
dimension \eqref{vd} to zero. The simplest way is by imposing incidence
conditions: we count only those curves which pass though fixed cycles
in $X$. 
Let 
$$
\mathcal C\to X\times\MM_g(X,\beta)
$$
be the universal curve.
We would like to intersect $\mathcal C$
 with a cycle $\alpha$ pulled back from $X$.  
Transversality issues again arise here,
so we use Poincar\'e dual cocyles.\footnote{Even if two submanifolds do not intersect transversally, the integral of the Poincar\'e dual cohomology class of one over the other still gives the correct topological intersection.} 
Let 
$$f\colon\mathcal C\to X \ \ \ \text{and} \ \ \ 
\pi\colon\mathcal C\to\MM_g(X,\beta)$$
be
the universal map and the 
projection to $\MM_g(X,\beta)$ respectively.
Let
$$
\widetilde\alpha=\pi_*\big(f^*\mathrm{PD}(\alpha)\big)\in H^*(\MM_g(X,\beta)).
$$
If $\alpha$ is a cycle of real codimension $a$, 
then $\widetilde\alpha$ is a cohomology class{\footnote{The
cohomological push-forward here uses the fact that
$\pi$ is an lci morphism. Alternatively flatness
can be used \cite{SGA}.}} in degree $a-2$ . 
When  transversality is satisfied, 
$\widetilde\alpha$
is Poincar\'e dual to the locus of curves in $\MM_g(X,\beta)$ which 
intersect $\alpha$.
After imposing sufficiently many incidence conditions to cut the
 virtual dimension to zero, we  define the {\em Gromov-Witten
invariant}
$$
N^{\scriptscriptstyle{\mathsf{GW}}}_{g,\beta}(\alpha_1,\ldots,\alpha_k)=
\int_{[\MM_g(X,\beta)]^{vir}}\widetilde\alpha_1
\wedge\ldots\wedge\widetilde\alpha_k\ \in \mathbb{Q}\,.
$$

We view the Gromov-Witten invariant{\footnote{We 
drop the superscript $\mathsf{GW}$ when clear from context.}}
 $N_{g,\beta}$ as
counting the  
curves in $X$ which pass through the cycles $\alpha_i$. 
The deformation invariance of $N_{g,\beta}$  follows
from construction of the virtual class. 
We are free to deform $X$ and the cycles $\alpha_i$ 
 in order to compute $N_{g,\beta}$.

The projective variety $X$ may be viewed as a symplectic manifold
with symplectic form obtained from the projective embedding. 
In fact, $N_{g,\beta}$ can be defined on any symplectic manifold $X$ by picking a compatible almost complex structure and using pseudo-holomorphic maps of curves. The resulting invariants do not depend on the choice of compatible almost complex structure, so define invariants of the symplectic
structure.\footnote{The role of the symplectic structure in the definition of the invariants is well hidden. Via Gromov's results,
the symplectic structure
is crucial for the compactness of the moduli space of stable maps.}

We can try to perturb
the  almost complex structures to make the moduli space transverse of 
the correct dimension. 
But even when embedded pseudo-holomorphic curves in $X$ are well-behaved,
their multiple covers invariably are not. Even within symplectic geometry,
the correct treatment of Gromov-Witten theory currently involves 
virtual classes.

% To get numerical invariants we need to cut down the virtual cycle from
% $vd$ dimension to zero. The traditional way to do this is to impose incidence
% conditions -- counting only those curves passing though some fixed cycle
% in $X$ -- by adding marked points to $C$ and insisting these lie on the cycle.
% So we introduce
% $$\M_{g,n}(X,\beta):=
% \left\{f\colon (C,p_1,\ldots,p_n)\to X\ \colon\ 
% \begin{array}{c}p_i\mathrm{\ distinct\ \emph{smooth\ points}} \\ \mathrm{of\ }C,\mathrm{\ and\ Aut}(f)
% \mathrm{\ finite} \end{array}\!\!\!\right\}\bigg/\cong\,.
% $$
% As before $C$ has arithmetic genus $g$, is at worst nodal (with none of
% the points $p_i$ being nodes), and Aut$(f)$ is the group of holomorphic
% automorphisms of $C$ which fix both $f$ and the marked points $p_i$.  Again
% this set has the structure of a Deligne-Mumford stack. It is compact; for
% instance if who points $p_i$ and $p_j$ move together on $C$ in $M_{g,n}(X,\beta)$
% then in the limit we attach a $\PP^1$ tail to $C$ at the limiting point
% and move the two marked points to distinct points on this $\PP^1$ (up to
% automorphisms of the $\PP^1$ it doesn't matter which ones). The resulting
% curve has 3 special points on the $\PP^1$ tail (the node where it attaches
% plus the two marked points) and so no extra automorphisms.

\subsection*{\S Advantages} 
Gromov-Witten theory is defined
for spaces $X$ of all dimensions and has been
proved to be 
a symplectic invariant 
(unlike most of the theories we will describe below).
As the first deformation invariant theory constructed, 
Gromov-Witten theory
 has been intensively studied for more than
20 years --- by now there are many exact
calculations and significant
structural results related to   
 integrable hierarchies and mirror symmetry. 

Since the moduli space of stable maps $\MM_g(X,\beta)$
lies over
the moduli space $\MM_g$ of stable
curves, Gromov-Witten theory is intertwined 
with the geometry of $\MM_g$.
Relations in the cohomology of $\MM_{g,n}$
yield universal differential equations for the
generating functions of Gromov-Witten invariants.
The most famous case is the WDVV equation \cite{DVV,WW} obtained
by the linear equivalence of the boundary strata of
$\MM_{0,4}$. 
The WDVV equation implies  the
associativity of the quantum cohomology ring of $X$
defined via the genus 0 Gromov-Witten invariants.
For example,
associativity for $\PP^3$ implies 
80160
twisted cubics  meet 12 general lines \cite{DI,FuPa}.
Higher genus relations such as
Getzler's \cite{GG} in genus 1  and the BP equation \cite{BPP}
in genus 2 also exist.

Gromov-Witten 
theory has links in many directions.
When $X$ is a curve, Gromov-Witten theory
is related to counts of Hurwitz covers \cite{OP2}.
For the Severi degrees of
curves in $\PP^2$ and $\PP^1 \times \PP^1$, Gromov-Witten
theory agrees with naive counts (when the
latter are sensible).
For surfaces of general type, 
Gromov-Witten theory links beautifully with 
Seiberg-Witten theory 
 \cite{Taubes}.
For 3-folds, there is a subtle and surprising relationship
between Gromov-Witten theory and the sheaf counting
theories discussed here in later sections.
The relation  with mirror symmetry \cite{CC,Giv,LLY} is a
high point of the subject.

\subsection*{\S Drawbacks} 
The theory is extremely hard to compute: 
even the Gromov-Witten theories of varieties of dimensions
0 and 1 are very complicated. The theory
of a point is related to the KdV hierarchy \cite{Wit}, and
the theory of $\PP^1$ is related to the
Toda hierarchy \cite{OP2}.
While such connections are 
beautiful,
using Gromov-Witten theory
to actually count curves is difficult, essentially due 
to the nonlinearity of maps from curves to varieties. 
The sheaf theories considered in the next sections 
concern more linear objects.

Because of the finite automorphisms of stable maps,
Gromov-Witten invariants are typically rational numbers. 
An old
idea in Gromov-Witten theory is that
underlying the rational Gromov-Witten invariants 
should be integer-valued curve counts. 
For instance, consider a stable map $f\in\MM_g(X,\beta)$ 
 double covering  
an image curve $C\subset X$ in class $\beta/2$. 
Suppose, for simplicity, $f$ and $C$ are rigid and unobstructed.
 Then, $f$ counts $1/2$ towards the Gromov-Witten 
invariant $N_{g,\beta}(X)$ because of its $\Z/2$-stabiliser.
Underlying this rational number is an integer $1$ counting the embedded curve $C$ in class $\beta/2$.

\subsection*{\S Serious difficulties} 
For the case of $3$-folds, 
Gromov-Witten theory is {\em not} enumerative
in the naive sense in genus $g>0$ 
due to degenerate contributions. 
The departure from naive counting happens already in positive genus 
for $\PP^3$.

Let $X$ be a $3$-fold.
The formula for the expected dimension of the moduli space
of stable maps \eqref{vd} is not genus dependent.
Consider a nonsingular embedded rigid rational curve 
\begin{equation} \label{grt23}
\PP^1\subset X
\end{equation} 
 in homology class $\beta$. 
The curve not only contributes $1$ to $N_{0,\beta}$,
 but also contributes in a 
complicated way to $N_{g\ge1,\beta}$. By attaching to 
the $\PP^1$ any stable curve $C$ at a nonsingular point, 
we obtain a stable map 
in the same class $\beta$ which collapses $C$ to a point. 
The 
contribution of \eqref{grt23}
to the Gromov-Witten invariants $N_{g\ge1,\beta}$ of $X$
 must be computed via 
integrals over the moduli spaces of stable curves. 
The latter integrals 
 are hard to motivate from the
point of view of curve counting.

A rather detailed study of the {\em Hodge integrals} over the
moduli spaces of curves
which arise in such degenerate contributions in Gromov-Witten theory
has been pursued  \cite{FaberPand,Pdegen}. A  main outcome has been an understanding
of the relationship of Gromov-Witten theory to naive curve counting
on $3$-folds in the Calabi-Yau  and Fano cases. 
The conclusion is a precise conjecture expressing integer
counts in terms of Gromov-Witten invariants (see the BPS conjecture
in the next section). 
The sheaf counting theories developed
later are now viewed as a more direct path to the integers
underlying Gromov-Witten theory in dimension 3.

What happens in higher dimensions? Results of \cite{KlemmPand,PandZing} for spaces
$X$ of dimensions
4 and 5 show a similar underlying integer structure for 
Gromov-Witten theory. However, a direct interpretation of the integer
counts (in terms of sheaves or other structures)
in dimensions higher than 3 awaits discovery.

\newsection{Gopakumar-Vafa / BPS invariants} \label{GVsec}
\subsection*{\S Invariants}
BPS invariants were introduced for Calabi-Yau 3-folds by Gopakumar-Vafa in \cite{GV1, GV2} using an M-theoretic construction.
The multiple
cover calculations \cite{FaberPand,Pdegen} in Gromov-Witten theory
provided basic motivation.
The definitions and conjectures related to
BPS states
were generalised to arbitrary 3-folds in \cite{Pdegen,PandICM}. 
While the original approach to the subject is not yet on 
a rigorous footing,
the hope is to define curve counting invariants  
 which avoid the
multiple cover and   degenerate contributions
of Gromov-Witten theory.
The BPS counts should  be the  integers underlying the 
rational Gromov-Witten invariants of $3$-folds.

To simplify the discussion here, let $X$ be a Calabi-Yau $3$-fold.
 Gopa\-kumar and Vafa consider a moduli space $\M$ of D-branes supported on curves in class $\beta$. While the precise mathematical definition is not clear, for a nonsingular embedded curve 
$C\subset X$ of genus $g$ and
class 
$$[C]=\beta\in H_2(X,\mathbb{Z})\,,$$ the D-branes are 
believed to be
(the pushforward to $X$ of) line bundles on $C$ of a
fixed degree, with moduli space a Jacobian torus 
diffeomorphic to $T^{2g}$. For singular curves,
the D-brane moduli space
should be a type of
relative compactified Jacobian over the ``space" of curves of class $\beta$.

Mathematicians have tended to interpret $\M$ as a moduli space of stable sheaves with 1-dimensional support in class $\beta$ and holomorphic Euler characteristic $\chi=1$. The latter condition 
is a technical device to rule out strictly
semistable sheaves. Over nonsingular curves $C$,
 the moduli space is simply  $\Pic_g(C)$. For singular curves,
more exotic sheaves in the compactified Picard scheme
arise. For nonreduced curves, we can find 
higher rank sheaves supported on the 
underlying reduced curve.
The support map
$\M\to B$, taking such a sheaf to the underlying 
support curve, is also required for the geometric path to 
the BPS invariants.
Here, $B$ is an appropriate (unspecified) parameter space of
curves in $X$.
For instance, there is certainly such a support 
map to the Chow variety of 1-cycles in $X$.

Let us now imagine that
we are in the ideal situation where the parameter space
$B=\coprod_iB_i$ is a disjoint union of connected components over which the map
 $\M\to B$ is a product, 
$$\M=\coprod\M_i\ \ \ \text{and}  \ \ \ \M_i=B_i\times F_i$$
with fibres $F_i$. 
The supposition is not ridiculous:
 the virtual dimension of curves in a Calabi-Yau 3-fold is 0,
 so we might hope that $B$ is a finite set of points. Then,
\beq{Kunneth}
H^*(\M)=\bigoplus_i H^*(B_i)\otimes H^*(F_i)\,.
\eeq
When each $B_i$ parameterises nonsingular curves of genus $g_i$ only, 
$$H^*(F_i)=H^*(T^{2g_i})=
(H^*(S^1))^{\otimes2g_i}$$ 
has normalised Poincar\'e polynomial 
$$P_y(F_i)=y^{-g_i}(1+y)^{2g_i}\,.$$
 Here, we normalise by shifting cohomological degrees by $-\dim_\C(F_i)$ to make $P_y(F_i)$ symmetric about degree 0. 
Then, $P_y(F_i)$ is a palindromic Laurent polynomial
 invariant under $y\leftrightarrow y^{-1}$ by Poincar\'e duality.

For more general $F_i$, the normalised Poincar\'e polynomial $P_y(F_i)$ is again invariant under $y\leftrightarrow y^{-1}$ if $H^*(F_i)$ satisfies even dimensional Poincar\'e duality.
Therefore, $P_y(F_i)$ may be written as a finite integral combination of terms $y^{-r}(1+y)^{2r}$, since the latter
 form a basis for the palindromic Laurent polynomials. 
Thus we can express $H^*(F_i)$ as a virtual combination of cohomologies of even dimensional tori. For instance, a cuspidal elliptic curve is topologically $S^2$ with 
$$P_y=y^{-1}(1+y^2)=(y^{-1}+2+y)-2=P_y(T^2)-2P_y(T^0)\,.$$
Cohomologically, we interpret the
cuspidal elliptic curve
 as 1 Jacobian of a genus 1 curve minus 2 Jacobians of genus 0 curves.

To tease the ``number of genus $r$ curves" in class $\beta$ from \eqref{Kunneth}, Gopakumar-Vafa write
\begin{equation}
\label{hqhq}
\sum_i(-1)^{\dim B_i}e(B_i)P_y(F_i) \quad\mathrm{as}\quad \sum_rn_r(\beta)y^{-r}(1+y)^{2r}
\end{equation}
and define the integers $n_r(\beta)$ to be the BPS invariants counting genus $r$ curves in class $\beta$. In Section \ref{MNOPsec},
 we will see that when $B$ is nonsingular
 and can be broken up into a finite number of points by a generic deformation, that number of points is $(-1)^{\dim B}e(B)$, see for instance \eqref{signed}. In other words, the virtual class of $B$ consists of $(-1)^{\dim B}e(B)$ points, explaining the first term in \eqref{hqhq}.

The K\"unneth decomposition \eqref{Kunneth} does {\em not} 
hold for general $\M\to B$, but can be replaced by the 
associated 
Leray spectral sequence. According to \cite{HST}, 
the \emph{perverse} Leray spectral sequence on intersection cohomology is preferable since it collapses and its terms satisfy the Hard Lefschetz theorem (which replaces the 
Poincar\'e duality used above). 
At least when $B$ is nonsingular,
 $\M$ is reduced with sufficiently mild singularities,
and $$\pi:\M\rightarrow B$$
 is equidimensional of fibre dimension $f=\dim\M-\dim B$,
 we can take
\beq{HST}
y^{-f}\sum_j(-1)^{\dim B}e({}^{p\!}R^j\pi_*\mathcal{IC}(\underline\C))y^j
=\sum_rn_r(\beta)y^{-r}(1+y)^{2r}
\eeq
as the Hosono-Saito-Takahashi definition\footnote{The original sources \cite{GV1, GV2, HST} make a great deal of use of $\mathfrak{sl}_2
\times\mathfrak{sl}_2$-actions on the cohomology of $\M$, but the end result is equivalent to the above intuitive description: decompose the fibrewise cohomology of $\M$ into the cohomologies of Jacobian tori, then take signed Euler characteristics in the base direction.}
of 
the BPS invariants
$n_r(\beta)$. 

The entire preceding discussion of BPS invariants
is only  motivational. We have not been precise about the
definition of the moduli space $B$. Moreover,
the hypotheses imposed in the
above constructions are rarely met (and when the
hypotheses fail, the constructions are usually
unreasonable or just wrong). Nevertheless,
there should exists BPS invariants
$n_{g,\beta}\in\Z$ ``counting'' curves of genus $g$
and class $\beta$ in $X$.

In addition to the $M$-theoretic construction, Gopakumar and
Vafa have made a beautiful prediction of
the relationship of the BPS counts to Gromov-Witten theory.
For Calabi-Yau $3$-folds, the conjectural formula is
\begin{multline}
\label{BPS}
\sum_{g\ge0,\ \beta\ne0}N_{g,\beta}^{\scriptscriptstyle {\mathsf{GW}}}
  u^{2g-2}v^\beta\ = \\
\sum_{g\ge0,\ \beta\ne0}n_{g,\beta}u^{2g-2}\sum_{d>0}\frac1d\left(
\frac{\sin(du/2)}{u/2}\right)^{2g-2}\!v^{d\beta}\,.
\end{multline}
The trigonometric terms on the right are motivated
by multiple cover formulas in Gromov-Witten theory \cite{FaberPand, Pdegen}.
The entire geometric discussion can be bypassed by {\em
defining} the BPS invariants via Gromov-Witten theory
by equation \eqref{BPS}. A
precise conjecture \cite{BryP1} then arises.

\vspace{15pt}
\noindent{\bf BPS conjecture I.} {\em For the  $n_{g,\beta}$
defined
via Gromov-Witten theory and formula \eqref{BPS},
the following properties hold:
\begin{enumerate}
\item[(i)] $n_{g,\beta}\in\Z$,
\item[(ii)]
for fixed
$\beta$, the $n_{g,\beta}$ vanish except for
finitely many $g\geq 0$.
\end{enumerate}}
\vspace{15pt}

For other $3$-folds $X$, when the virtual dimension is
positive
$$\int_\beta c_1(X)>0\,,$$
incidence conditions to cut down the virtual dimension 
to 0
must be included. This case will be
discussed in Section \ref{Kim} below.
The conjectural formula for the BPS counts
is similar, see \eqref{BPS+}.

\subsection*{\S Advantages} For 3-folds, BPS invariants should be the ideal curve counts. The BPS invariants are integer valued and coincide
with naive counts in many cases where the latter make sense.
For example, the BPS counts (defined via Gromov-Witten theory)
agree with naive curve counting in $\PP^3$ in genus 0, 1, and 2.
The definition via Gromov-Witten theory shows 
$n_{g,\beta}$ is a  symplectic invariant.

For Calabi-Yau 3-folds $X$, the  BPS counts do not always agree
with naive counting.
A trivial example is the slightly different treatment of
an embedded super-rigid elliptic curve $E\subset X$, see \cite{Pdegen}.
Such an $E$ contributes a single BPS count to each  multiple
degree $n[E]$. A much more subtle BPS contribution
is given by a super-rigid genus 2 curve C in class $2[C]$
\cite{BryP2,BryP3}. 
We view BPS counting now as more fundamental than naive
curve counting (and equivalent to, but {\em not} always equal to,
naive counting).

\subsection*{\S Drawbacks} 
The main drawback is the murky foundation of the
geometric construction of the BPS invariants.
For nonreduced curves, the contributions of the higher rank
moduli spaces of sheaves on the underlying support
curves  remain mysterious.
The real strength of the theory will only be realised
after the foundations are clarified. For example,
properties (i) and (ii) of the BPS conjecture should
be immediate from a geometric construction.
The definition via Gromov-Witten theory is far from
adequate.

A significant limitation of the BPS counts is
the restriction to 
3-folds. However,
calculations \cite{KlemmPand, PandZing} show some hope
of parallel structures in
higher dimensions, see also
\cite{JoyceK}.

\subsection*{\S Serious difficulties}
The geometric foundations appear very hard to establish.
There is no likely path in sight (except in genus 0 where Katz has made a 
rigorous proposal \cite{KatzBPS}, see Section \ref{PTsec}).
The Hosono-Saito-Takahashi approach does not incorporate
the virtual class (the term $(-1)^{\dim B}e(B)$ is a crude approximation for
 the virtual class of the base $B$) and fails in
general.

%Using a map between moduli spaces, and how this behaves with respect to cohomology ($\mathfrak sl_2\times\mathfrak sl_2$-actions, or the perverse filtration) seems a very long way from our current methods of constructing invariants. In particular the HST construction did not incoorporate a virtual cycle  and so is not yet deformation invariant, contradicting \eqref{BPS}.

Developments concerning 
motivic invariants \cite{JoySong, KS} and the  categorification of  
invariants with cohomology theories instead of Euler characteristics
 appear somewhat closer to the methods required in the Calabi-Yau 3-fold case.
 For instance, Behrend has been working to categorify 
his constructible function \cite{BehrendDT} to give a 
perverse sheaf that could replace $\mathcal{IC}
(\underline\C)$ in the HST definition, perhaps yielding
a deformation invariant theory. Even then, why formula
\eqref{BPS} should hold is a mystery.

An approach to BPS invariants via  stable
pairs (instead of Gromov-Witten theory) will be discussed
in Section \ref{PTsec} below. The BPS invariants $n_{g,\beta}$ are
there
again defined  
by a formula similar to \eqref{BPS}. 
 The
stable pairs perspective is better than the
Gromov-Witten approach and has led to substantial 
 recent progress 
 \cite{ChDiacPan, MaulikYun, MigShende, Shende, TodaBPS}.
Nevertheless, the hole in the subject left by the lack of
a direct geometric construction is not yet filled.

\newsection{Donaldson-Thomas theory} \label{MNOPsec}
\subsection*{\S Moduli}
Instead of considering maps of curves into $X$, we
 can instead study embedded curves.
Let a {\em subcurve} $Z\subset X$ be a subscheme of
dimension 1.
The \emph{Hilbert scheme} 
 compactifies embedded curves 
by allowing them to degenerate to arbitrary subschemes. 
Let 
$I_n(X,\beta)$ be the Hilbert scheme parameterising
subcurves $Z\subset X$  with
$$\chi(\O_Z)=n \in \mathbb{Z} \ \ \ \text{and} \ \ \ [Z]=\beta\in H_2(X)\,. $$
Here, $\chi$ denotes the holomorphic
Euler characteristic and $[Z]$ denotes the class of the
subcurve (involving only the 1-dimensional components). 
By the above conditions, $I_n(X,\beta)$ parameterises subschemes
which are unions of possibly nonreduced curves and points in $X$.

We give a few examples to show how the
Hilbert scheme
 differs from the space of stable maps. First,
 consider a family of nonsingular conics
\beq{conics}
C_{t\ne0}=\{x^2+ty=0\}\subset\C^2
\eeq
as a local model which can, of course,
 be further embedded in any higher dimension. The natural limit as $t\to0$,
\begin{equation}\label{n55}
C_0=\{x^2=0\}\subset\C^2,
\end{equation}
is indeed the limit in the Hilbert scheme.
\begin{figure}
\vspace{3mm}
\input{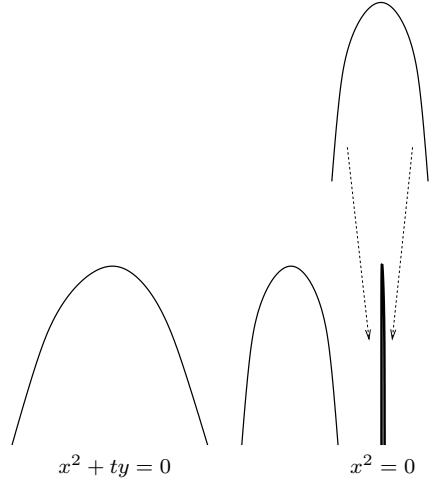}
\caption{The degeneration \eqref{conics} with the limiting stable map double covering $x=0$.}
\end{figure}
The limit \eqref{n55}
 is the $y$-axis with multiplicity two thickened
 in the $x$-direction. 

In the stable map case,
the limit of the family \eqref{conics} is very different. 
There we take the limit of the associated map from 
$\C$ to $C_t$ given by\footnote{The formula gives a well-defined map only modulo automorphisms of the curve --- specifically the automorphism $\xi\mapsto-\xi$.}
$$\xi\mapsto (-t^{1/2}\xi,\xi^2)\,.$$
 The result is the \emph{double cover} $\xi\mapsto(0,\xi^2)$ of the $y$-axis. 
So the thickened scheme in the Hilbert scheme
is replaced by the double cover. The latter
is an orbifold point in the space of stable maps with $\Z/2$-stabiliser given by $\xi \mapsto -\xi$. \medskip
\begin{figure}
\input{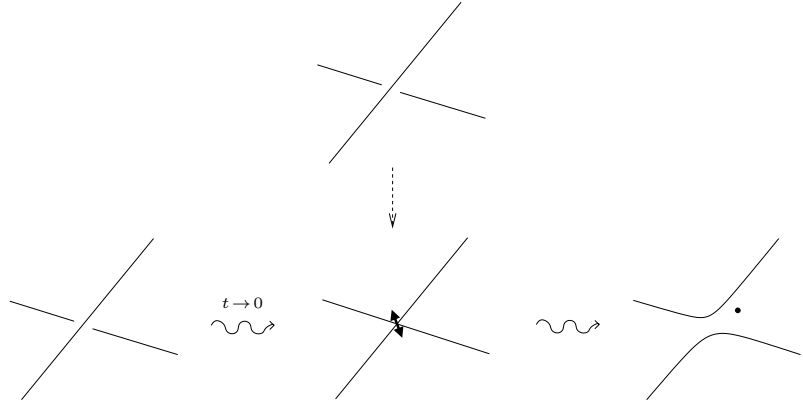}
\caption{The family \eqref{genuschange} with the
subscheme limit below and the stable map limit above. 
On the right is a deformation of the limit subscheme with a free point breaking off.} \label{gc}
\end{figure}

In the next example, we illustrate the phenomenon of genus 
change which occurs only in dimension at least 3. 
A global model is given by a twisted rational cubic in $\PP^3$
degenerating 
to a plane cubic of genus 1 \cite{PS}. An easier local model $C_t\subset\C^3$ has 2 components: the $x$-axis in the plane $z=0$, and the $y$-axis moved up into the plane $z=t$,
\beq{genuschange}
C_t=\{x=0=z\}\sqcup\{y=0=z-t\}\subset\C^3\,,
\eeq
see Figure \ref{gc}. As a stable map, we take the associated
 inclusion of two copies of the affine line $\C$. The stable
map limit at $t=0$ takes
the same domain $\C\sqcup\C$ onto the $x$- and $y$- axes, an embedding away from the origin where the map is 2:1. In other words, the limit stable
map\footnote{There is another stable map given by the embedding of the image
\eqref{r45}. In a compact global model,
 the latter would be a map from a nodal stable curve of genus one larger
so would \emph{not} feature in the compactification of the family we are considering.} is the normalisation of the image 
\begin{equation}\label{r45}
\{xy=0=z\}\subset\C^3\,.
\end{equation}

In the Hilbert scheme, the limit of the
family \eqref{genuschange}
is rather worse. The ideal of $C_t$ is
$$
(x,z)\cdot(y,z-t)\ =\ \big(xy,x(z-t),yz,z(z-t)\big)\,.
$$
We take the limit as $t\to0$. The flat limit
here happens to be the ideal generated by the limit of the above generators. 
The limit ideal does not contain $z$:
\beq{embedded}
(xy,xz,yz,z^2)\ \subsetneq\ (xy,z) \,.
\eeq
However, 
 after multiplying $z$ by any element of the maximal ideal $(x,y,z)$ 
of the origin, we land inside the limit ideal. 
Therefore, the limit curve  is given by $\{xy=0=z\}$
with a scheme-theoretic embedded point added at the origin pointing along the $z$-axis --- in the direction along which the two components came together. 
The embedded point ``makes up for" the point lost in the intersection
and ensures that the family of curves is flat over $t=0$.

In a further flat family, the embedded
 point can  break off, and the curve can be smoothed $\{xy=\epsilon,\,z=0\}$ to a curve of higher genus. 
In the Hilbert scheme, 
we have all 1-dimensional subschemes made up of curves \emph{and} points, 
with curves of different genus balanced by extra free points. 
The constant $n$ in $I_n(X,\beta)$ is $1-g+k$ for a reduced curve of arithmetic genus $g$ with $k$ free and embedded points added, 
so we can increase $g$ at the expense of increasing $k$ by the same amount.

\subsection*{\S Deformation theory}
Hilbert schemes of curves can have 
arbitrary dimensional components and terrible singularities.
Worse still, 
the natural deformation/obstruction theory of the Hilbert scheme
does not lead to a virtual class. 
However, if we restrict attention to 3-folds $X$  and view
$I_n(X,\beta)$ as a moduli space of {\em sheaves}, then we obtain 
a different obstruction theory which does admit a virtual class.
The latter observation is the starting point of Donaldson-Thomas
theory.

Given a 1-dimensional subscheme $Z\subset X$, the associated
ideal sheaf $\I_Z$ is a stable sheaf 
with Chern character $$(1, 0,-\beta,-n)\in
H^0\oplus H^2\oplus H^4\oplus H^6$$ and trivial determinant. 
Conversely, all such stable sheaves with trivial determinant 
can be shown to embed in their double duals $\O_X$ and
thus are all ideal sheaves. Hence,
$I_n(X,\beta)$ is a moduli space of sheaves, at least set theoretically.
With more work, an isomorphism of schemes
can be established, see \cite[Theorem 2.7]{PT1}.

The moduli space of sheaves $I_n(X,\beta)$ also admits a virtual class
\cite{ThCasson, MNOP1}.
The main point is that deformations and obstructions are governed by
\beq{exts}
\Ext^1(\I_Z,\I_Z)_0 \quad\mathrm{and}\quad \Ext^2(\I_Z,\I_Z)_0
\eeq
respectively, where the subscript $0$ denotes the trace-free part governing deformations with \emph{fixed determinant}. 
Since $\Hom(\I_Z,\I_Z)=\C$
consists of only the scalars, the 
trace-free part vanishes. 
By Serre duality,
\begin{equation*}
\Ext^3(\I_Z,\I_Z)\cong\Hom(\I_Z,\I_Z\otimes K_X)^*\cong H^0(K_X)^* \cong H^3(\O_X)\,.
\end{equation*}
The last group $H^3(\O_X)$
is removed when taking trace-free parts. 
Hence, the terms \eqref{exts} are the only nonvanishing trace-free Exts, 
and there are {\em no} higher obstruction spaces. 
The Exts \eqref{exts} 
 govern a perfect obstruction theory of virtual dimension equal to
$$
\text{ext}^1(\I_Z,\I_Z)_0-\text{ext}^2(\I_Z,\I_Z)_0=\int_\beta c_1(X)\,,
$$
compare \eqref{vd}. If the virtual dimension
is positive, insertions are needed  
to produce invariants \cite{MNOP1}.

On Calabi-Yau 3-folds, 
moduli of sheaves admit a particularly 
nice deformation-obstruction theory \cite{DT, ThCasson}.
The deformation and obstruction spaces \eqref{exts} are \emph{dual} to each other,
\beq{duality}
\Ext^2(\I_Z,\I_Z)_0\cong\Ext^1(\I_Z,\I_Z)_0^*\,,
\eeq
by Serre duality.
Any moduli space of sheaves on a
 Calabi-Yau 3-fold can be realized as
the critical locus of a holomorphic function on an
ambient nonsingular space: the holomorphic Chern-Simons functional 
in infinite dimensions \cite{W,DT} or locally on an appropriate 
finite dimensional slice \cite{JoySong}. Since the moduli space 
is the zero locus of a closed 1-form, the obstruction space is 
the cotangent space at any point of moduli space. 
More generally, Behrend \cite{BehrendDT} calls obstruction theories 
satisfying the global version of \eqref{duality} {\em symmetric}.
The condition
is equivalent to asking for the moduli space to be locally the zeros of 
an {\em almost closed} 1-form on a smooth ambient space ---
 a 1-form with 
exterior derivative vanishing scheme theoretically on the moduli space.\footnote{By \cite{1form},
the condition is strictly weaker than asking for the moduli space to be locally the zeros of a closed 1-form.}

If the moduli space of sheaves is nonsingular
 (but of too high dimension), then the symmetric obstruction theory 
forces the obstruction bundle to be globally isomorphic to the cotangent bundle. The virtual class, here the top Chern class of the obstruction bundle,
is then just the signed topological Euler characteristic of the moduli space
\beq{signed}
(-1)^{\dim I_n(X,\beta)}e(I_{n}(X,\beta))\,.
\eeq
Remarkably, 
Behrend shows that for any moduli space $\M$ with 
a symmetric obstruction theory there is a constructible function 
$$\chi^B\colon \M\to\Z$$ 
with respect to which the weighted Euler characteristic gives the integral
of the virtual class \cite{BehrendDT}. 
Therefore, each point of the moduli space contributes in a local way to the global invariant, by $(-1)^{\dim \M}$ for a nonsingular point and by
a complicated number taking multiplicities into account for singular points. 
When $\M$ is locally the critical locus of a function, the number 
is $(-1)^{\dim\M}(1-e(F))$ where $F$ is the Milnor fibre of our point.
Unfortunately, how to find a parallel approach to the virtual class
when $X$ is not Calabi-Yau is not currently known.

Integration against  the virtual class of
$I_n(X,\beta)$ yields the
Donald\-son-Thomas invariants. In the Calabi-Yau case, no
insertions are required:
$$I_{n,\beta}\ =\ \int_{[I_n(X,\beta)]^{vir}} 1\ =\ 
e\big(I_n(X,\beta),\chi^B\big)\,. $$
Since $I_{n}(X,\beta)$ is a scheme (ideal sheaves
have no automorphisms) and $[I_n(X,\beta)]^{vir}$
is a cycle class with $\mathbb{Z}$-coefficients, the invariants
$I_{n,\beta}$ are {\em integers}.
Deformation invariance of $I_{n,\beta}$ follows from
properties of the virtual class.

\subsection*{\S MNOP conjectures}
A series of conjectures linking the Donaldson-Thomas theory of
3-folds to Gromov-Witten theory were advanced in \cite{MNOP1,MNOP2}.
For simplicity, we restrict ourselves here to the Calabi-Yau
case.

For fixed curve class $\beta\in H_2(X,\mathbb{Z})$,
the Donaldson-Thomas {\em partition function} is 
$$
\mathsf{Z}^{\mathrm{DT}}_\beta(q)=\sum_nI_{n,\beta}q^n\,.
$$
Since $I_n(X,\beta)$ is easily seen to be empty for
$n$ sufficiently negative, the partition function is
a Laurent series in $q$. 
To count just curves, and not points and curves, 
MNOP form the \emph{reduced} generating function \cite{MNOP1} by 
dividing by the contribution of just points:
\beq{reduced}
\mathsf{Z}_\beta^{\mathrm{red}}(q)=
\frac{\mathsf{Z}_\beta^{\mathrm{DT}}(q)}{\mathsf{Z}^{\mathrm{DT}}_0(q)}\,.
\eeq
MNOP first conjectured the degree $\beta=0$ contribution can be calculated as
$$
\mathsf{Z}^{\mathrm{DT}}_0(q)=M(-q)^{e(X)}\,,
$$
where $M$ is the MacMahon function,
$$
M(q)=\prod_{n\ge1}(1-q^n)^{-n}\,,
$$
the generating function for 3d partitions. 
Proofs can now be found in  \cite{BFHilb,LiDT,LPCobordism}. 
Second, MNOP conjectured $\mathsf{Z}_\beta^{\mathrm{red}}(q)$ is 
the Laurent expansion of a \emph{rational function} in 
$q$, invariant\footnote{The Laurent series itself need not be $q\leftrightarrow q^{-1}$ invariant. For instance the rational function $q(1+q)^{-2}$ is invariant, but the associated Laurent series $q-2q^2+3q^3-\ldots$ is not.} under $q\leftrightarrow q^{-1}$.
Therefore, we can substitute $q=-e^{iu}$ and obtain
 a real-valued function of $u$. The main conjecture of
MNOP in the Calabi-Yau case is the following.

\vspace{15pt}
\noindent{\bf GW/DT Conjecture:} 
$\mathsf{Z}^{\mathrm{GW}}_\beta(u)=
\mathsf{Z}^{\mathrm{red}}_\beta(-e^{iu})$.
\vspace{15pt}

The conjecture 
asserts a precise equivalence relating Gromov-Witten to 
Donaldson-Thomas theory.
Here,
 $$\mathsf{Z}^{\mathrm{GW}}_\beta(u)=
\sum_{g\ge0}N_{g,\beta}^{\bullet} \ u^{2g-2}$$ 
is the generating function of \emph{disconnected} Gromov-Witten 
invariants $N_{g,\beta}^\bullet$ defined just as in Section \ref{GW} by relaxing the condition that the curves be connected, 
but excluding maps which contract connected components to points. 
Equivalently, $\mathsf{Z}^{\mathrm{GW}}_\beta(u)$ is
 the exponential of the generating function of 
connected Gromov-Witten invariants $N_{g,\beta}$,
$$
\sum_{\beta\ne0} \mathsf{Z}^{\mathrm{GW}}_\beta(u)v^\beta=\sum_{\beta\ne0,\,g\ge0}
N_{g,\beta}^\bullet \ u^{2g-2}v^\beta=
\exp\bigg(\sum_{\beta\ne0,\,g\ge0}N_{g,\beta}\,u^{2g-2}v^\beta\bigg).
$$
A version of the GW/DT correspondence with insertions for 
non Calabi-Yau 3-folds can be found in \cite{MNOP2}. Various refinements, involving theories relative to a divisor, or equivariant with respect to a group action, are also conjectured.
All of these conjectures  
have been proved for toric 3-folds in \cite{MOOP}.

The GW/DT conjecture 
should be viewed as
involving an analytic continuation and series
expansion about two different points ($q=0$ and $q=-1$, corresponding to $u=0$). Therefore, the conjecture
cannot be
understood term by
term\footnote{When combined with the Gopakumar-Vafa formula \eqref{BPS}
 and the relationship to the stable pairs discussed below, 
the GW/DT conjecture
 will become rather more comprehensible, see \eqref{bpsform}.} ---
  to determine a single invariant on one side of the conjecture, 
knowledge of all of the invariants on the other side is necessary.

The overall shape of the conjecture is clear:
the two different ways of counting curves in a fixed class $\beta$ are entirely equivalent, with \emph{integers} determining the Gromov-Witten invariants of 3-folds. By \cite[Theorem 3.20]{PT1}, the integrality prediction of the GW/DT correspondence is entirely equivalent to the integrality prediction of the Gopakumar-Vafa formula \eqref{BPS}.

\subsection*{\S Advantages} The integrality of the invariants
is a significant advantage of 
using the Hilbert schemes $I_n(X,\beta)$ to define a
counting theory. Also,
the virtual counting of subschemes, at least in the Calabi-Yau
3-fold case,
fits into the larger context of
 counting higher rank bundles, sheaves, and objects of the derived category
of $X$. The many recent developments in wall-crossing \cite{JoySong, KS} 
apply to this more general setting. 
We will see an example in the next section.

Behrend's constructible function sometimes
makes computations (in the Calabi-Yau case at least) more feasible ---  we can use cut and paste techniques to reduce to more local calculations. See for instance \cite{BehrendBryan}.

\subsection*{\S Drawbacks} The theory only works for nonsingular projective varieties of dimension at most  3.
While the Hilbert scheme of curves is  always well-defined,
the deformation/obstruction theory fails
to be 2-term in higher dimensions.
By contrast, Gromov-Witten theory is well-defined
in all dimensions and is proved to be a symplectic
invariant. While we expect Donaldson-Thomas theory to have
a fully symplectic approach,
how to proceed is not known. 

In Gromov-Witten theory, the genus expansion 
makes a connection to the
moduli of curves (independent of $X$).
The Euler characteristic $n$ plays a parallel role
in Donaldson-Thomas theory, but is much less useful.
While there are very good low genus results
in Gromov-Witten theory, there are few  analogues for
the Hilbert scheme.

Behrend's constructible function approach for
the Calabi-Yau case is difficult to use. For example,
the constructible functions even for
toric Calabi-Yau 3-folds have not been determined.{\footnote{
Amazingly, we do not even know whether the
constructible functions are non-constant in the
toric Calabi-Yau case!}}
So far,  Behrend's theory has been useful mainly for
formal properties related to motivic invariants
and wall-crossing. For more concrete calculations
 involving Behrend's functions see \cite{BehrendBryan, BBS}.

\subsection*{\S Serious difficulties} 
For the GW/DT correspondence,
the division by $\mathsf{Z}_0(q)$ confuses the
geometric interpretation of the invariants.
In fact, the subschemes of $X$ with free points 
make the theory rather unpleasant
to work with. This ``compactification" of the space of embedded curves 
is much larger than the original space, adding enormous
components with free points. In practice, the free points lead to
constant technical headaches (which play little role in
the main development of the invariants).

It is tempting to think of working with the closure of 
the ``good components" of the Hilbert scheme instead, but such 
an approach would not have a reasonable
 deformation theory nor a virtual class. 
However, a certain birational modification of the idea does work and
will be discussed in the next section.

% more linear but no longer 1 dim so virt cycle problems in dims $>3$.
% Now relative theory \cite{LiRelative}. 

\newsection{Stable pairs} \label{PTsec}
\subsection*{\S Limits revisited} 
Consider again the family of Figure  \ref{gc}.
For $t\neq0$,
denote the disjoint union \eqref{genuschange} by 
$$C_t=C^1_t\cup C^2_t\,.$$
The ideal sheaf $\I_{C_t}$,
central to the Hilbert scheme analysis,  is just the kernel of the surjection
\beq{pairt}
\O_X\to\O_{C^1_t}\oplus\O_{C^2_t}\,.
\eeq
%The data \eqref{pairt} is equivalent to the ideal sheaf. 
For the moduli of stable pairs,
the map itself (not just the kernel) will be 
fundamental. 
We will take a 
natural limit of the map given by 
\beq{pair}
\O_X\to\O_{C_0^1}\oplus\O_{C_0^2}\ 
\eeq
where the limits of the component curves are
$$C^1_0=\{x=0=z\} \ \ \ \text{and} \ \ \  C^2_0=\{y=0=z\}\,.$$
 The result is a map which is \emph{not} a surjection at the origin (where $C_1$ and $C_2$ intersect and the sheaf on the right has rank 2). 
In the
 limit, there is a nonzero cokernel, the structure sheaf of the origin $\O_0$, which accounts for the extra point lost in the intersection. 
Losing surjectivity
replaces the embedded point arising in the limit of ideal sheaves \eqref{embedded}.

The cokernels of the above maps \eqref{pairt} are \emph{not} flat over $t=0$
even though the sheaves
$\O_{C^1_t}\oplus\O_{C^2_t}$ \emph{are} flat. 
Similarly the kernels of the maps \eqref{pairt} are \emph{not} flat over $t=0$.
In fact, at $t=0$, we get the ideal $(xy,z)$ of $C^1_0\cup C^2_0$ 
which we already saw in \eqref{embedded} is not the flat limit of the ideal sheaves of $C_t$.

\subsection*{\S Moduli} 
The limit \eqref{pair} is an example of a stable pair.
The moduli of stable pairs provides a different
 sheaf-theoretic compactification of the space of embedded curves. 
The moduli space is intimately related to the Hilbert scheme, 
but is much more efficient.

Let $X$ be a nonsingular projective $3$-fold.
A {\em{stable pair}} $(F,s)$ is a coherent sheaf $F$ with dimension 1 support in $X$ and a  section $s\in H^0(X,F)$ satisfying the following stability condition:
\begin{itemize}
\item $F$ is \emph{pure}, and
\item the section $s$ has zero dimensional cokernel.
\end{itemize}
Let $C$ be the scheme-theoretic support of $F$.
Condition (i) means all the irreducible components
of $C$ are of dimension 1 (no 0-dimensional components).
By  \cite[Lemma 1.6]{PT1}, $C$ has no embedded points.
A stable pair $$\O_X\to F$$ 
therefore
defines a Cohen-Macaulay curve $C$ via the kernel $\I_C\subset \O_X$
 and a 0-dimensional subscheme of $C$ via the support of the
 cokernel\footnote{When $C$ is Gorenstein (for instance if $C$ lies
in a nonsingular surface), stable pairs supported on $C$ are in bijection with 0-dimensional subschemes of $C$. More precise scheme theoretic isomorphisms of moduli spaces are proved in \cite[Appendix B]{PT3}.}.

To a stable pair, we associate the Euler characteristic and
the class of the support $C$ of $F$,
$$\chi(F)=n\in \mathbb{Z} \  \ \ \text{and} \ \ \ [C]=\beta\in H_2(X,\mathbb{Z})\,.$$
For fixed $n$ and $\beta$,
there is a projective moduli space of stable pairs $P_n(X,\beta)$ \cite[Lemma 1.3]{PT1} by work of Le Potier \cite{LePot}. 
While the Hilbert scheme $I_n(X,\beta)$ is a moduli space of curves plus free and embedded points, $P_n(X,\beta)$
 should be thought of as a moduli space of curves plus points \emph{on the curve} only. 
Even though points
still play a role (as the example \eqref{genuschange} shows),
the moduli of stable pairs is
much smaller than $I_n(X,\beta)$.

\subsection*{\S Deformation theory}
To define a flexible counting theory, 
a compactification of the family of curves in $X$
should admit a 2-term deformation/obstruction theory and a virtual class. 
As in the case of $I_n(X,\beta)$, the most immediate
 obstruction theory of $P_n(X,\beta)$ does \textit{not} admit such a structure. For $I_n(X,\beta)$, a
solution was found by considering a subscheme $C$ to be equivalent to a sheaf $\I_C$ with trivial determinant. For $P_n(X,\beta)$, we consider a stable pair to define an object of $D^b(X)$, the quasi-isomorphism equivalence class of the complex
\beq{Idot}
I\udot=\{\O_X\Rt{s}F\}\,.
\eeq
For $X$ of dimension 3, the object $I\udot$ determines the stable pair \cite[Proposition 1.21]{PT1}, and the fixed-determinant deformations of $I\udot$ in $D^b(X)$ match those of the pair $(F,s)$ to all orders
\cite[Theorem 2.7]{PT1}. The latter property shows 
the scheme $P_n(X,\beta)$ may be viewed as a moduli space
of objects in the derived 
category.{\footnote{Studying the moduli of objects
in the derived category is a young subject. Usually, such constructions
lead to Artin stacks. The space $P_n(X,\beta)$ is a
rare example where a component of the moduli of objects
in the derived category is a scheme (uniformly for all $3$-folds $X$).}}
We can then use the obstruction theory of the complex
$I\udot$ in place of the
obstruction theory of the pair.

The deformation/obstruction theory
for complexes  is governed at  $[I\udot]\in P_n(X,\beta)$ by
\beq{exts2}
\Ext^1(I\udot,I\udot)_0 \quad\mathrm{and}\quad \Ext^2(I\udot,I\udot)_0\,.
\eeq
Formally, the outcome is parallel
to \eqref{exts}.
The obstruction theory \eqref{exts2} has all the attractive properties
of the Hilbert scheme case: 2 terms, 
a virtual class of dimension $\int_\beta c_1(X)$, and a description 
via the  $\chi^B$-weighted Euler characteristics in the Calabi-Yau case.

\subsection*{\S Invariants} After imposing incidence conditions (when
$\int_\beta c_1(X)$ is positive)
and integrating against the virtual class, we obtain
stable pairs invariants for $3$-folds $X$.
In the Calabi-Yau case, the invariant is  the length of the virtual class:
$$
P_{n,\beta}=\,\int_{[P_n(X,\beta)]^{vir}}1\ =\ e\big(P_n(X,\beta),\chi^B\big)\,.
$$
For fixed curve class $\beta\in H_2(X,\mathbb{Z})$,
the stable pairs {partition function} is 
$$
{\mathsf{Z}}^{\mathrm{P}}_{\beta}(q)=\sum_nP_{n,\beta}q^n.
$$
Again, elementary arguments show the moduli spaces
$P_n(X,\beta)$ are empty for sufficiently negative $n$, so
${\mathsf{Z}}^{\mathrm{P}}_{\beta}$ is a Laurent series in $q$.
Since the free points are now confined to the curve instead of 
roaming over $X$, we do not have to form a reduced series as in 
\eqref{reduced}. In fact, we conjecture \cite[Conjecture 3.3]{PT1} 
the partition function
$\mathsf{Z}^{\mathrm{P}}_\beta$ to be precisely the reduced theory of Section \ref{MNOPsec}.

\vspace{15pt}
\noindent{\bf DT/Pairs Conjecture:} 
$
{\mathsf{Z}}^{\mathrm{P}}_{\beta}(q)\ =\ {\mathsf{Z}}^{\mathrm{red}}_{\beta}(q).
$
\vspace{15pt}

The DT/Pairs correspondence is expected for all 3-folds $X$
with the incidence conditions playing no significant role \cite{MNOP2}.
Using the definition 
$\mathsf{Z}^{\mathrm{red}}_{\beta}=\mathsf{Z}^{\mathrm{DT}}_\beta/
\mathsf{Z}_0^{\mathrm{DT}}$,
we find
\beq{wallnos}
\sum_m   P_{n-m,\beta}\cdot    I_{m,0}  =I_{n,\beta}\,.
\eeq
Relation \eqref{wallnos} should be interpreted as a wall-crossing formula for
counting invariants in the derived category of coherent sheaves $D^b(X)$
under a change of stability condition. 

For invariants of Calabi-Yau 3-folds,
wall-crossing has been studied intensively in recent years, 
  and we give only the briefest of descriptions.
Ideal sheaves parameterised by $I_n(X,\beta)$ are Gieseker stable. 
We can imagine changing the stability condition\footnote{Ideally, we would 
work with Bridgeland stability conditions \cite{Br}, but that 
 is not currently possible. The above discussion can be made precise using the limiting stability conditions of \cite{BayerPoly,TodaStab}, or even Geometric Invariant Theory \cite{StopTh}.} to destabilise the ideal shaves 
with free and embedded points.
If $Z$ is a 1-dimensional subscheme, then $Z$
 has a maximal pure dimension 1 subscheme $C$ defining a sequence
$$
0\to\I_Z\to\I_C\to Q\to0\,,
$$
where $Q$ is the maximal 0-dimensional subsheaf of $\O_Z$. In $D^b(X)$,
 we  equivalently have  the exact triangle
\beq{destab}
Q[-1]\to\I_Z\to\I_C\,.
\eeq
We can imagine  the stability condition crossing a wall 
on which the phase (or slope) of $Q[-1]$ equals that of $\I_C$. 
On the other side of the wall, $\I_Z$ will be destabilised 
by \eqref{destab}. 
Meanwhile,  
extensions $E$ in the opposite direction
\beq{restab}
\I_C\to E\to Q[-1]
\eeq
will become stable.
But stable pairs are just such extensions! The exact sequence 
$$0\to\I_C\to\O_X
\rt{s}F\to Q\to0$$ yields the exact triangle
$$
\I_C\to I\udot\to Q[-1]\,.
$$
The
moduli space of pairs $P_{n}(X,\beta)$ should 
give precisely the space of stable objects for the new 
stability condition. 

The formula \eqref{wallnos} for $I_{n,\beta}-P_{n,\beta}$ should follow from the more general wall-crossing formulae 
of \cite{JoySong, KS}. The $m^{th}$ term in (\ref{wallnos}) is the correction from subschemes
$Z$ whose maximal 0-dimensional subscheme (or total number of free
and embedded points) is of length $m$. It involves both the space $\PP(\Ext^1(\I_C,Q[-1]))$ of extensions \eqref{destab} and the space $\PP(\Ext^1(Q[-1],\I_C))$ of extensions \eqref{restab}. Though both are hard to control, they contribute to the wall-crossing formula through the difference in their Euler characteristics\footnote{Really, we need to weight by the restriction of the Behrend function $\chi^B$. To make the above analysis work then requires $\chi^B$ to satisfy the identities of \cite{JoySong, KS}. In fact, the automorphisms of $Q$ make the matter much more complicated than we have suggested.}, which is the topological number 
$$\chi(\I_C,Q)=\ \,{\text{length}}\,(Q)\ =m\,.$$ 

The above sketch has now been carried out at the level of (unweighted) Euler characteristics \cite{TodaDTPT, StopTh} and for
$\chi^B$-weighted Euler characteristics in \cite{BridgelandDTPT} in the Calabi-Yau case. 
The upshot is the DT/Pairs conjecture is now proved for Calabi-Yau 3-folds.
 The rationality of 
$\mathsf{Z}_\beta^{\mathrm{red}}(q)$ 
and the symmetry  under $q\leftrightarrow q^{-1}$ 
is also proved \cite{BridgelandDTPT}.

\subsection*{\S Example}
Via the Behrend weighted Euler characteristic approach to 
 the invariants of a Calabi-Yau 3-fold, we
can talk about the contribution of a single curve $C\subset X$ to 
the stable pairs generating function $\mathsf{Z}^P_\beta(q)$.
No such discussion is possible in Gromov-Witten theory.

If $C$ is nonsingular of genus $g$, then the stable pairs 
supported on $C$ with $\chi=1-g+n$ are in bijection with $\Sym^nC$ 
via the map taking a stable pair to the support of the cokernel $Q$.
Therefore, $C$ contributes\footnote{The Behrend function restricted to $\Sym^nC$ can be shown \cite[Lemma 3.4]{PT3} to be the constant $(-1)^{n-g}c$, where $c=\chi^B(\O_C)$ is the Behrend function of the moduli space of torsion sheaves evaluated at $\O_C$.}
\beq{sym}
\mathsf{Z}^{\mathrm{P}}_C(q)\,=\,c\sum_n(-1)^{n-g}e(\Sym^nC)q^{1-g+n}\,=\,(-1)^gc\,q^{1-g}(1+q)^{2g-2}.
\eeq
The rational function on the right is
invariant under $q\leftrightarrow q^{-1}$. We view the
symmetry as a manifestation of Serre duality (discussed below).
Control of the free points in stable pair theory makes the geometry more transparent. The same calculation for $\mathsf{Z}^{\mathrm{red}}_C(q)$ 
based on the Hilbert scheme is much less enlightening.
The above calculation is closely related
to the BPS conjecture for stable pairs.

\subsection*{\S Stable pairs and BPS invariants} 

By a formal argument \cite[Section 3.4]{PT1},  
the stable pairs partition function can be written uniquely
in the following special way:
\begin{align*}
\!\!{\mathsf{Z}}^{\mathrm{P}}(q,&v):=1+\sum_{\beta\ne0}
\mathsf{Z}^{\mathrm{P}}_{\beta}(q)v^\beta \\ & = \exp\bigg(\sum_{r}
\sum_{\gamma\neq 0}\sum_{d\geq 1} \tilde{n}_{r,\gamma}\frac{(-1)^{(1-r)}}d
\ (-q)^{d(1-r)}(1-(-q)^d)^{2r-2} v^{d\gamma}\bigg)\,,\!\!\! 
\end{align*}
where the $\tilde{n}_{r,\gamma}$ are integers and
vanish for fixed $\gamma$ and $r$ sufficiently large.

We can compose the various conjectures to link
the BPS counts of Gopakumar and Vafa to the stable pairs invariants. 
The form we get  from the conjectures is almost exactly as above:
\begin{multline*}
{\mathsf{Z}}^{\mathrm{P}}(q,v)
  = \\ \exp\bigg(\sum_{r\geq 0}
\sum_{\gamma\neq 0}\sum_{d\geq 1} n_{r,\gamma}\frac{(-1)^{(1-r)}}d
\ (-q)^{d(1-r)}(1-(-q)^d)^{2r-2} v^{d\gamma}\bigg)\,.\!\!\!
\end{multline*}
The only difference is the restriction on the $r$ summation.
Hence, we can {\em define} the BPS state counts by stable pairs
invariants via the $\tilde{n}_{r,\gamma}$!

\vspace{15pt}
\noindent{\bf BPS conjecture II.} {\em For the  $\tilde{n}_{r,\beta}$
defined
via stable pairs theory, the vanishing
$$\tilde{n}_{r<0,\beta} = 0$$
holds for $r<0$.}
\vspace{15pt}

By its construction, the approach to defining the  BPS states counts 
via stable pairs satisfies the full integrality condition
and half of the finiteness of BPS conjecture I.
We therefore regard the stable pairs
perspective as better than the path via Gromov-Witten
theory. Still, a direct construction of the BPS invariants
along the lines discussed in Section \ref{GVsec} would be
best of all.\footnote{For curves with only reduced plane curve singularities, both constructions of BPS numbers have been shown to coincide after making the $\chi^B=(-1)^{\dim}$ approximation to the virtual class \cite{MaulikYun, MigShende}.}

\vspace{10pt}
For irreducible classes\footnote{A class $\beta\in H_2(X,\mathbb{Z})$ 
is irreducible if it cannot be written as a sum $\alpha+\gamma$ of nonzero classes containing algebraic
 curves.},
the BPS formula for the stable pairs invariants  can be written as
\beq{bpsform}
\mathsf{Z}^{\mathrm{P}}_{\beta}(q)=\sum_{r\geq 0}^gn_{r,\beta}\,q^{1-r}(1+q)^{2r-2},
\eeq
with $n_{r,\beta}=0$ for all sufficiently large $r$.
There is a beautiful interpretation of
\eqref{bpsform}
in the light of \eqref{sym}: to the stable pairs invariants, the curves in class $\beta$ look like a disjoint union of a finite number $n_{r,\beta}$ of 
nonsingular curves of genus $r$.

We can prove directly that the  partition function  
$\mathsf{Z}_\beta^{\mathrm{P}}$ can be written in the form \eqref{bpsform}.
% We do so at the level of Euler characteristics for simplicity; incorporating
% the Behrend function is not much harder \cite{PT3}.
For $r\ge1$, the functions $q^{1-r}(1+q)^{2r-2}$,
$$
1,\qquad q^{-1}+2+q,\qquad q^{-2}+4q^{-1}+6+4q+q^2,\qquad q^{-3}+\ldots
$$
form a natural $\Z$-basis for the Laurent polynomials invariant under $q\leftrightarrow q^{-1}$. For $r=0$, the coefficients of the Laurent series do not satisfy the same symmetry,
$$
q(1+q)^{-2}\ =\ q-2q^2+3q^3-4q^4+\ldots\ =\ \sum_{n\ge1}(-1)^{n-1}nq^n.
$$
To prove \eqref{bpsform}, it is therefore equivalent to  show  the coefficients $P_{n,\beta}$ of the partition
function satisfy not the $q\leftrightarrow q^{-1}$ symmetry but
\beq{Serre}
P_{n,\beta}=P_{-n,\beta}+c(-1)^{n-1}n
\eeq
for some constant $c$.

Relation \eqref{Serre} is
 a simple consequence of Serre duality for the fibres of the Abel-Jacobi map. By forgetting the section,
 we obtain a map from stable pairs to stable sheaves\footnote{The irreducibility of $\beta$ implies the sheaves with arise 
 are stable since $F$ has rank 1 on its irreducible support.},
$$
\begin{array}{ccc}
P_n(X,\beta) & \To & \M_n(X,\beta)\,,\vspace{7pt} \\
(F,s) & \mapsto & F\,.
\end{array}
$$
The fibre of the map is $\PP(H^0(F))$ with weighted Euler characteristic\footnote{As proved in  \cite[Theorem 4]{PT3}, the Behrend function is constant on $\PP(H^0(F))$
 with value $(-1)^{n-1}c$ where $c=\chi^B(\O_C)$. On a first reading,
the Behrend function can be safely ignored here.} $(-1)^{n-1}c\cdot
h^0(F)$. 
%We
%write $F=\iota_*L$ where 
%$$\iota\colon C=\mathrm{Supp}(F)\into X$$
% and $L$ is a rank one sheaf on $C$. Then 
There is an isomorphism
$$
\begin{array}{ccc}
\M_n(X,\beta) & \To & \M_{-n}(X,\beta)\,, \vspace{7pt}\\
F &\mapsto& F^\vee\,.
\end{array}
$$
where $F^\vee=\ext^2(F,K_X)$. If $F$ is the push-forward
of a line bundle $L$ from a nonsingular curve $C$, then 
$F^\vee$ is the push-forward of $L^*\otimes \omega_C$, see \cite{PT3}
for details.
The fibre $\PP(H^0(F^\vee))$ over $F^\vee$ is $\PP(H^1(F)^*)$ by Serre duality,
 with weighted Euler characteristic $(-1)^{-n-1}c\cdot h^1(F)$.

To prove relation \eqref{Serre}, we calculate the
difference between the two above contributions:
$$
(-1)^{n-1}c(h^0(F)-h^1(F))=(-1)^{n-1}c\chi(F)=(-1)^{n-1}cn\,.
$$
Summation over the space of stable sheaves (in the sense of Euler characteristics) yields the relation
\beq{SD}
P_{n,\beta}-P_{-n,\beta}=(-1)^{n-1}n\,e(\M_n(X,\beta),c)\,.
\eeq
The weighted Euler characteristics
\beq{TodaConj}
e(\M_n(X,\beta),c)=e(\M_{n+1}(X,\beta),c)
\eeq
are independent of $n$:  tensoring with a 
degree 1 line bundle 
 relates sheaves supported on $C$ with $\chi=n$ to those with $\chi=n+1$. 
We have proved the relation \eqref{Serre}. \medskip

The above argument 
shows the coefficient $n_{0,\beta}$ of $q(1+q)^{-2}$ is the $\chi^B$-weighted Euler characteristic of $\M_n(X,\beta)$.
In fact, Katz \cite{KatzBPS} had previously proposed 
the DT invariant of $\M_1(X,\beta)$ as a good definition of $n_{0,\beta}$ 
for \emph{any} class $\beta$, not necessarily irreducible. 
Naively, Katz's definition sees only the rational curves
 because for a curve of higher genus the action of the
 Jacobian on the moduli space of sheaves forces the 
(weighted) Euler characteristic of the latter to be zero. 
Katz's proposal can be viewed as a weak analogue of the genus by genus
methods in Gromov-Witten theory. 
\medskip

Identity \eqref{SD} is easily seen to be another wall-crossing formula \cite{BayerPoly, PT3, TodaStab}. In \cite{TodaSurvey}, Toda has extended the above analysis to all curve classes by 
extending the methods of Joyce \cite{JoySong} and the ideas of Kontse\-vich-Soibelman \cite{KS} on BPS formulations of general sheaf counting. His main
result reduces  BPS conjecture II 
to an analogue of identity \eqref{TodaConj} for DT invariants 
for dimension 1 sheaves for \emph{all} classes $\beta$.\footnote{When $\beta$ is not irreducible, sheaf stability issues change the definition of the DT invariant, see \cite[Conjecture 6.3]{TodaSurvey} for details.}

\subsection*{\S Advantages} The stable pair theory has the advantages of the ideal sheaf theory --  integer invariants conjecturally equivalent to the rational Gromov-Witten invariants -- but with the bonus of eliminating the free points on $X$. The geometry of the BPS conjectures is more clearly explained
by stable pairs than any of the other approaches. 

If descendent insertions (coming from higher Chern classes of 
the tautological bundles) are considered, the theory of stable
pairs behaves much better than the parallel constructions
for the moduli of stable maps or the Hilbert scheme.
For example, the descendent partition functions for stable
pairs are rational in $q$. See \cite{PaPix} for proofs in toric cases
and further discussion.

At least for 3-folds, stable pairs appears to be the best
counting theory to consider at the moment. The main hope
for a better approach lies in  the
direct geometric construction of the BPS counts.

\subsection*{\S Drawbacks} Just as for the
Donaldson-Thomas theory of ideal sheav\-es, 
the stable pairs invariants
 have only been constructed on nonsingular projective varieties of 
dimension at most 3. While we expect a
parallel theory for symplectic invariants, 
how to proceed is not clear.

\subsection*{\S Serious difficulties} 
In the theory of stable pairs,
free points are allowed to
move along the support curve $C$. 
The free points are
necessary to probe the geometry of the curve (and the
associated BPS contributions in all genera) but in a rather roundabout way. 
An alternative opened up by the Behrend function might be to work with 
open moduli spaces (on which the arithmetic genus does not jump), but deformation invariance then becomes problematic.

\begin{figure}
\vspace{5mm}
\input{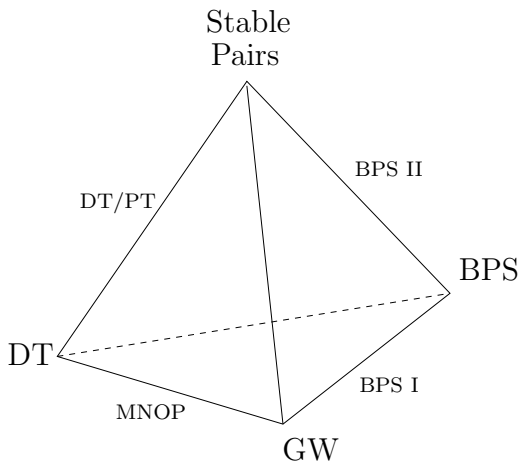}
\caption{Conjectures connecting curve counting theories}
\end{figure}

As we have said repeatedly,
a rigorous and sheaf theoretic approach to BPS invariants (at
least for 3-folds and possibly in higher dimensions as well)  
 would be highly desirable.

\newsection{Stable unramified maps} \label{Kim}

\subsection*{\S Singularities of maps}
A difficulty which arises in Gromov-Witten theory
is the abundance of collapsed components.
In the moduli space of higher genus stable maps to $\PP^1$ of
degree $1$, the entire complexity
comes from such collapsed components attached to
a degree 1 map  of a genus 0 curve to $\PP^1$. 
Collapsed contributions have to be removed
to arrive at the integer counts underlying Gromov-Witten theory.

A map $f$ from a nodal curve $C$ to a nonsingular variety $X$ is
{\em unramified} at a nonsingular point  $p\in C$ if the differential 
$$df: T_{C,p} \rightarrow T_{X,f(p)}$$
is injective. If $f$ is unramified at $p$, the component
of $C$ on which $p$ lies cannot be collapsed.

The idea of stable unramified maps, introduced 
by  Kim,  Kresch, and  Oh \cite{KKO}, is to control both the
domain (allowing only nodal curves) {\em and}  the
singularities of the maps (essentially unramified 
and with no collapsed components). 
The price for these properties is paid in the
complexity of the target space $X$. The target
cannot remain inert, but must be allowed to
degenerate.

\subsection*{\S Degenerations}
Let $X$ be a nonsingular projective variety of dimension $n$.
The Fulton-MacPherson \cite{FuMa} configuration space 
$X[k]$ compactifies the moduli of $k$ distinct labelled
points on $X$. The Fulton-MacPher\-son compactification  may be
viewed as a higher dimensional analogue of the geometry of
marked points on stable curves --- when the points attempt
to collide, the space $X$ degenerates and the colliding
points are separated in a bubble. 

The possible degenerations of $X$ which occur are
easy to describe. We start with the trivial family
$$\pi:X\times \bigtriangleup_0  \rightarrow \bigtriangleup_0\,,$$
with fibre $X$ over the disk $\bigtriangleup_0$  with base point $0$.
Next, we allow an iterated sequence of finitely many 
blow-ups of the total
space $X\times \bigtriangleup_0$ at points which, at each stage,
\begin{enumerate}
\item[(i)] lie over $0\in \bigtriangleup_0$ and
\item[(ii)] lie in the smooth locus of the morphism to $\bigtriangleup_0$.
\end{enumerate}
After the sequence of blow-ups is complete, we take the
fibre $\widetilde{X}$ of the resulting total space 
$\widetilde{X\times \bigtriangleup_0}$
over $0\in \bigtriangleup_0$.
The space $\widetilde{X}$, a {\em Fulton-MacPherson degeneration} of
$X$, is a normal crossings divisor in the total space.

The Fulton-MacPherson degeneration $\widetilde{X}$
 contains a distinguished
component $X_+$ which is a blow-up of the original $X$ at distinct points.
The other components of $\widetilde{X}$ are simply blow-ups of
$\PP^n$. Of the latter, there are two special types
\begin{enumerate}
\item[(i)] {\em ruled} components ($\PP^n$ blown-up at 1 point),
\item[(ii)] {\em end} components ($\PP^n$ blown-up at 0 points).
\end{enumerate}
The singularities of $\widetilde{X}$ occur only in the
intersections of the components.

By construction, there is a canonical morphism 
$$\rho:\widetilde{X} \rightarrow X$$
which blows-down $X_+$ and contracts the other components
of $\widetilde{X}$.
The automorphisms of $\widetilde{X}$ which commute with $\rho$
can only be non-trivial on the components of type (i) and (ii).

For the moduli of stable unramified maps, the target $X$ is
allowed to degenerate to any Fulton-MacPherson degeneration
$\widetilde{X}$.

\subsection*{\S Moduli}
Let $X$ be a  nonsingular projective variety of dimension $n$.
The moduli space $\M_g(X,\beta)$ of {\em stable unramified maps}
to $X$ parameterises the data
$$C \stackrel{f}{\rightarrow} \widetilde{X} \stackrel{\rho}{\rightarrow}
X$$
satisfying the following conditions:
\begin{enumerate}
\item[(i)] $C$ is a connected nodal curve of arithmetic genus $g$,
\item[(ii)] $\widetilde{X}$ is a Fulton-MacPherson degeneration of
$X$ with canonical contraction $\rho$,
\item[(iii)] $\rho_*f_*[C] = \beta \in H_2(X,\mathbb{Z})$,
\item[(iv)] the nonsingular locus of $\widetilde{X}$ pulls-back
            to exactly the nonsingular locus of $C$,
$$f^{-1}(\widetilde{X}^{ns}) = C^{ns}\,,$$
\item[(v)] $f$ is {\em unramified} on $C^{ns}$,

\item[(vi)] at each node $q\in C$, the two incident
branches $B_1,B_2 \subset C$ map to two {\em different} components 
$Y_1,Y_2 \subset \widetilde{X}$ and meet the  
intersection divisor at $q$ with equal multiplicities,
$$\Big[ B_1 \cdot Y_1 \cap Y_2\Big]_{Y_1,q} =
\Big[ B_2 \cdot Y_1 \cap Y_2\Big]_{Y_2,q}\,, $$
\item[(vii)]
for each ruled component $R\subset \widetilde{X}$, there is
a component of $C$ which is mapped by $f$ to $R$ with image
{\em not} equal to a fibre of the ruling,
\item[(viii)] for each end component $E\subset \widetilde{X}$,
there is
a component of $C$ which is mapped by $f$ to $E$ with image
{\em not} equal to a straight line.
\end{enumerate}

By (v), the map $f$ is unramified everywhere except
possibly at the nodes of $C$ (which must map to the
singular locus of $\widetilde{X}$).
Constraint (vi) is the standard {\em admissibility} condition
for infinitesimal smoothing which arises in relative
Gromov-Witten theory \cite{IPark, LRuan, lirel}.
Conditions (vii) and (viii) serve to stabilize the
components of $\widetilde{X}$ with automorphisms over $\rho$.

The moduli space $\M_{g}(X,\beta)$ of unramified maps
is a proper Deligne-Mumford stack. The unramified map limits of our two simple examples of degenerations (\ref{conics}) and (\ref{genuschange}) 
are easily described. 
For \eqref{conics}, the stable map limit is a double cover which is ramified
over two branch points. In the unramified limit,
 we take the Fulton-MacPherson degeneration which blows up these points 
in $X$ and adds projective space components. The proper transform of 
the double cover is then attached to nonsingular
 plane conics in the two added projective spaces. The conics are
tangent to the intersection divisors at the points hit by the double cover.
For \eqref{genuschange}, the limit is the same as in Gromov-Witten theory: the normalisation of the image \eqref{r45} in the trivial Fulton-MacPherson degeneration of $X$.

A central result of \cite{KKO} is the identification of
the deform\-ation/obstruc\-tion theory of an unramified
map mixing the (unobstruct\-ed) deformation theory of
Fulton-MacPherson degenerations with the usual
deformation theory of maps to $X$. The deformation/obstruction
theory is 2-term, and a virtual class is constructed 
on $\M_{g}(X,\beta)$ of dimension
$$\int_\beta c_1(X)+(\dim_\com X-3)(1-g) $$
as in Gromov-Witten theory.

There is no difficulty to include marked points in the
definition of unramified maps \cite{KKO}. Via incidence
conditions imposed at the markings, a full set of
unramified invariants can be constructed for any $X$.

\subsection*{\S Connections to BPS counts: CY case}
How do the unramified
invariants relate to all the other counting theories we have
discussed? Since unramified invariants have been introduced
very recently, not many calculations have been done.
In the case of Calabi-Yau 3-folds $X$, an attempt \cite{Set}
at finding
the analogue of the Aspinwall-Morrison formula for
multiple covers of an embedded $\PP^1\subset X$
with normal bundle $\mathcal{O}_{\PP^1}(-1) \oplus
\mathcal{O}_{\PP^1}(-1)$
showed the invariant was different for double covers.

A full transformation relating the unramified theory to
the other Calabi-Yau counts has not yet been proposed. 
Surely such a transformation exists
and has an interesting form.

\vspace{+15pt}
\noindent{\bf Question:}
What is the relationship between unramified invariants and
Gromov-Witten theory for the Calabi-Yau 3-folds?
\vspace{+15pt}

\subsection*{\S Connections to BPS counts: positive case}
Let $X$ be a nonsingular projective 3-fold and let
$\beta \in H_2(X,\mathbb{Z})$ be a curve class
satisfying
\begin{equation}\label{g553}
\int_\beta c_1(X) >0\,.
\end{equation}
Let 
$\gamma_1, \ldots, \gamma_n \in H^*(X,\mathbb{Z})$
be integral cohomology classes Poincar\'e dual to 
cycles in $X$ defining incidence conditions for curves.
We require
the dimension constraint
$$n+ \int_\beta c_1(X) = \sum_{i=1}^n \text{codim}_\com(\gamma_i)$$
to be satisfied.
Let 
$$N_{g,\beta}^{\scriptscriptstyle\mathrm{UR}}(\gamma_1,\ldots,\gamma_n) \in \mathbb{Q}$$ 
be the corresponding genus $g$ 
unramified invariant.

%\vspace{+15pt}
%\noindent{\bf Integrality conjecture:}
%$N_{g,\beta}^{unram}(\gamma_1,\ldots,\gamma_n) \in \mathbb{Z}\,.$
%\vspace{+15pt}

The BPS state counts
of Gopakumar and Vafa were generalized from the Calabi-Yau
to the positive case in \cite{PandICM, Pdegen}.
The BPS invariants 
$n_{g,\beta}(\gamma_1,\ldots, \gamma_n)$ are {\em defined}
via Gromov-Witten theory by:
\begin{align} \nonumber
\sum_{g\geq 0}
N_{g,\beta}^{\scriptscriptstyle\mathrm{GW}}(\gamma_1,&\ldots,\gamma_n)
\ u^{2g-2}= \\ \label{BPS+}
& \sum_{g\geq 0}
 n_{g,\beta}(\gamma_{1}, \cdots,
\gamma_{n})\
 u^{2g-2} \left( \frac{\sin
({u/2})}{u/2}\right)^{2g-2+\int_{\beta}c_1(X)}.
\end{align}
Zinger \cite{Zing} proved the above
definition yields {\em integers}
$n_{g,\beta}(\gamma_1,\ldots, \gamma_n)$
which vanish for sufficiently high $g$ (depending upon $\beta$)
when the positivity \eqref{g553} is satisfied.
The following conjecture{\footnote{BPS
conjecture III for unramified invariants was made by R.P.
and appears in Section 5.2 of \cite{KKO}.}}
connects the unramified theory to BPS counts.

\vspace{+15pt}
\noindent{\bf BPS conjecture III:}
$N_{g,\beta}^{\scriptscriptstyle\mathrm{UR}}(\gamma_1,\ldots,\gamma_n) =
n_{g,\beta}(\gamma_1,\ldots, \gamma_n)$.
\vspace{+15pt}

The above simple BPS relation should be true
because the moduli space of unramified maps
avoids all collapsed contributions. If proved,
unramified maps may be viewed as providing
a direct construction of the BPS counts
in the positive case.

\subsection*{\S Advantages} 
The main advantage of the unramified theory is the
simple form of the singularities of the maps. In particular,
avoiding collapsed components leads to (the expectations of)
much better behavior than Gromov-Witten theory.

The theory also enjoys many of the advantages of Gromov-Witten theory:
definition in all dimensions, relationship to the moduli of
curves, and connection with naive enumerative geometry for $\PP^2$ 
and $\PP^1 \times \PP^1$.

\subsection*{\S Drawbacks} 
The Fulton-MacPherson degenerations add a great deal
of complexity to calculations in the unramified theory.
Even in modest geometries, a large number of components
in the degenerations are necessary.
In localization formulas, Hodge integrals
on various Hurwitz/admissible cover moduli spaces occur
(analogous to the standard Hodge integrals on the moduli space
of curves appearing in Gromov-Witten theory).
While the latter have been studied for a long time,
the structure of the former has not been so carefully
understood.

Unramified maps remove the degenerate contributions of Gromov-Witten theory, but keep the multiple covers. For Calabi-Yau 3-folds,
 the invariants are rational numbers. 
The BPS invariants are expected to
underlie the theory, but how is not yet understood.

The unramified theory is expected to be symplectic, but
the details have not been worked out yet.

\subsection*{\S Serious difficulties} 
The theory has been studied for only a short time.
Whether the complexity of the degenerating
target is too difficult to handle remains to be seen.

\newsection{Stable quotients} \label{MOP}

\subsection*{\S Sheaves on curves}
We have seen compactifications of the family of curves on
$X$ via maps of nodal curves to $X$ and via sheaves on $X$.
The counting theory obtained from the
moduli space of stable quotients \cite{MOP}, involving
sheaves on nodal curves, takes a hybrid approach.
The stable quotients invariants are
directly connected to Gromov-Witten theory in many
basic cases. 
 However, the main application of  stable
quotients to date has been to the geometry of
the moduli space of curves.

\subsection*{\S Moduli} 
Let $(C,p_1,\ldots,p_n)$ be a connected nodal curve with nonsingular marked points.
Let $q$ be a quotient of the rank $N$ trivial bundle
 $C$,
\begin{equation*}
\com^N \otimes \oh_C \stackrel{q}{\rarr} Q \rarr 0\,.
\end{equation*}
If the quotient sheaf $Q$ 
is locally free at the nodes  of $C$,
 then
$q$ is a {\em quasi-stable quotient}. 
 Quasi-stability of $q$ implies the associated
kernel,
\begin{equation*}
0 \rightarrow S \rightarrow
\com^N \otimes \oh_C \stackrel{q}{\rarr} Q \rarr 0\,,
\end{equation*}
is a locally free sheaf on $C$. Let $r$ 
denote the rank of $S$.

Let $C$ be a curve equipped
with a quasi-stable quotient $q$.
The data $(C,q)$ determine 
a {\em stable quotient} if
the $\mathbb{Q}$-line bundle 
\begin{equation}\label{aam}
\omega_C(p_1+\ldots+p_n)
\otimes (\wedge^{r} S^*)^{\otimes \epsilon}
\end{equation}
is ample 
on $C$ for every strictly positive $\epsilon\in \mathbb{Q}$.
Quotient stability implies
$2g-2+n \geq 0$.

Viewed in concrete terms, no amount of positivity of
$S^*$ can stabilize a genus 0 component 
$$\proj^1\stackrel{\sim}{=}P \subset C$$
unless $P$ contains at least 2 nodes or markings.
If $P$ contains exactly 2 nodes or markings,
then $S^*$ {\em must} have positive degree.

\subsection*{\S Isomorphism} 
Two quasi-stable quotients on a fixed curve $C$
\begin{equation}\label{fpp22}
\com^N \otimes \oh_C \stackrel{q}{\rarr} Q \rarr 0,\qquad
\com^N \otimes \oh_C \stackrel{q'}{\rarr} Q' \rarr 0
\end{equation}
are {\em strongly isomorphic} if
the associated kernels 
$$S,S'\subset \com^N \otimes \oh_C$$
are equal.

An {\em isomorphism} of quasi-stable quotients
 $$\phi:(C,q) \rarr
(C',q)
$$ is
an isomorphism of curves
$$\phi: C \stackrel{\sim}{\rarr} C'$$
with respect to which 
the quotients $q$ and $\phi^*(q')$ 
are strongly isomorphic.
Quasi-stable quotients \eqref{fpp22} on the same
curve $C$
may be isomorphic without being strongly isomorphic.

The moduli space of stable quotients 
$\overline{\QQ}_{g}({\mathbb{G}}(r,N),d)$ parameterising the
data
$$(C,\  0\rarr S \rarr
\com^N\otimes \oh_C \stackrel{q}{\rarr} Q \rarr 0)\,,$$
with {\em rank}$(S)=r$ and {\em deg}$(S)=-d$,
is a proper Deligne-Mumford stack of finite type
over $\com$. A proof, by Quot scheme methods,
is given in \cite{MOP}.

Every stable quotient $(C,q)$
yields a rational map from the underlying curve
$C$ to the Grassmannian $\mathbb{G}(r,N)$.
If the quotient sheaf $Q$ is locally free on all of $C$, then
the stable quotient yields a regular map
from $C$ to the Grassmannian. Hence, we may view
stable quotients as compactifying the space of maps
of genus $g$ curves to Grassmannians of class $d$ times 
a line.

\subsection*{\S Deformation theory}
The moduli of stable
quotients 
 maps 
to the Artin stack of pointed domain curves
$$\nu^A:
\overline{\QQ}_{g}({\mathbb{G}}(r,N),d) \rightarrow {\mathfrak{M}}_{g,n}\,.$$
The moduli  of stable quotients with fixed underlying
curve 
$$[C] \in {\mathfrak{M}}_{g,n}$$
 is simply
an open set of the Quot scheme of $C$. 
The deformation theory of the Quot scheme 
determines a 2-term obstruction theory on
$\overline{\QQ}_{g}({\mathbb{G}}(r,N),d)$ relative to
$\nu^A$
given by $(R\Hom(S,Q))^\vee$.

More concretely, for the stable quotient,
\begin{equation*}
0 \rightarrow S \rightarrow
\com^N \otimes \oh_C \stackrel{q}{\rarr} Q \rarr 0\,,
\end{equation*} the
deformation and obstruction spaces relative to $\nu^A$ are
$\text{Hom}(S,Q)$ and $\text{Ext}^1(S,Q)$
respectively. Since $S$ is locally free and $C$ is a curve, the higher obstructions
$$\text{Ext}^{k}(S,Q)= H^{k}(C,S^*\otimes Q) = 0, \ \ \ k>1$$
vanish.

A quick calculation shows the virtual dimension of the moduli
of stable quotients 
equals the virtual dimension of the moduli of
stable maps to ${\mathbb{G}}(r,N)$.

\subsection*{\S Invariants.}
There is no difficultly in adding marked points
to the moduli of stable quotients, see \cite{MOP}.
Therefore, we can define a theory of stable quotients 
invariants for Grassmannians.
Similar targets such as flag varieties
for $\mathbf{SL}_n$ admit a parallel development.
An enumerative theory of stable quotients was sketched in
\cite{MOP} for complete intersections in such spaces. Hence,
there is a stable quotients theory for the Calabi-Yau
quintic in $\PP^4$.

Since \cite{MOP}, the construction of
stable quotient invariants has been extended to toric
varieties \cite{CK} and appropriate GIT quotients \cite{CKM}.
The associated counting theories (well-defined with
2-term deformation/ob\-struct\-ion theories) should
be regarded as depending not only on the target
space, but {\em also} on the quotient presentation.
The direction is related to 
the young subject of
gauged Gromov-Witten theory (and in particular to
the rapidly developing study of theories of
Landau-Ginzburg type \cite{ChangLi, FAN}).

\vspace{+15pt}
\noindent{\bf Question:}
What is the relationship between stable quotient invariants and
Gromov-Witten theory for varieties?
\vspace{+15pt}

For all flag varieties for $\mathbf{SL}_n$, the above question
has a simple answer: the counting by stable quotients
and Gromov-Witten theory agree exactly \cite{MOP}.
Perhaps exact agreement also holds for Fano toric varieties,
see the conjectures in \cite{CK}. But in the non-Fano cases,
and certainly for the Calabi-Yau quintic, the stable
quotient theory is very different. There should be
a wall-crossing understanding \cite{TodaSQ} of the transformations,
but much work remains to be done.

\subsection*{\S Advantages}
Stable quotients provide a more efficient compactification than
Gromov-Witten theory. In the case of projective space,
there is a blow-down morphism
$${\MM}_g(\PP^{N-1},d) \rightarrow \overline{\QQ}_g(
{\mathbb{G}}(1,N),d)$$
which pushes-forward the virtual class of the moduli of
stable curves to the virtual class of the moduli of
stable quotients \cite{MOP}.
A principal use of the moduli of stable quotients
has been to explore the tautological rings of the moduli
of curves \cite{PaPix2} --- and in particular to prove the Faber-Zagier
conjecture for relations 
among the $\kappa$ classes on $\mathcal{M}_g$ \cite{PaPa}.
The efficiency of the boundary
plays a crucial role in the analysis. 

The difference between stable maps and stable quotients
can be seen already for elliptic curves in projective
space. For stable maps,
the associated moduli space is singular with
multiple components. A desingularization, by 
{\em blowing-up},
is described in \cite{vak-zing} and applied to calculate
the genus 1 Gromov-Witten invariants of the quintic
Calabi-Yau in \cite{Zing2}. On the other hand, the moduli of stable quotients
related to such elliptic curves is a nonsingular {\em blow-down}
of the stable maps space \cite{MOP}. The stable quotients moduli here
is a much smaller compactification.{\footnote{A geometric
investigation by Cooper of the stable quotients spaces
in genus 1 for projective spaces can be found in \cite{cooper}.}}  A parallel application to
the genus 1 stable quotients invariants of the quintic
Calabi-Yau is a very natural direction to pursue.

\subsection*{\S Drawbacks}
The stable quotients approach to the enumeration of curves,
while valid for different dimensions, appears
to require more structure on $X$ (embedding, toric, or
quotient presentations). The method is therefore not
as flexible as Gromov-Witten theory.

Also, unlike Gromov-Witten theory, there is not
yet a symplectic development. However, the
connections with gauged Gromov-Witten theory may
soon provide a fully symplectic path to 
stable quotients.

\begin{figure}
\vspace{5mm}
\input{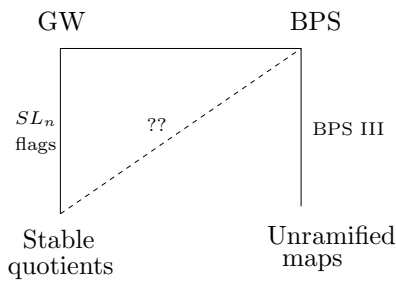}
\caption{Conjectures relating curve counting theories}
\end{figure}

\subsection*{\S Serious difficulties} 
The theory has been studied for only a short time.
The real obstacles, beyond those discussed
above, remain to be encountered.

\section*{Appendix: Virtual classes}
\setcounter{section}{7}

\subsection*{\S Physical motivation}
There are countless ways
to compactify the  spaces of curves in a projective variety $X$. 
What distinguishes the 6 main approaches we have described is the
presence in each case of a 
virtual fundamental class.

Moduli spaces arising in physics should naturally 
carry virtual classes when
cut out by a section (the derivative of an action functional) of a vector 
bundle over a nonsingular ambient space (the space of fields).
While both the space and bundle are usually infinite dimensional, the derivative of the section is often Fredholm, so 
we can make sense of the difference in the dimensions.
 The difference is  
the \emph{virtual dimension} of the moduli space
 -- the number of equations minus the number of unknowns. The question, though, of what geometric objects
to place in the boundary is often not so clearly specified in
the physical theory.

As an example, the space of $C^\infty$-maps from a Riemann surface $C$ to $X$, modulo diffeomorphisms of $C$, is naturally an infinite dimensional orbifold away from the maps with infinite automorphisms. Taking $\dbar$ of such a map gives a Fredholm section of the infinite rank bundle with fibre $\Gamma(\Omega^{0,1}_C(f^*T_X))$ over the map $f$. The zeros of the section
 are the holomorphic maps $$f\colon C\to X.$$ 
However,
to arrive at the definition of a stable
map requires further insights about nodal curves.

The 
Fredholm property allows us to take slices to reduce locally to the following finite dimensional model of the moduli problem.

\subsection*{\S Basic model}
Consider a nonsingular ambient variety $A$ of dimension $n$.
Let  $E$ be a rank $r$ bundle on $A$ with section $s\in\Gamma(E)$ with zero locus  $\M$:
\beq{model}
\spreaddiagramcolumns{-2pc}
\xymatrix{
& E\dto \\
\M=Z(s)\ \subset & A\ar@/^{-1.5ex}/[u]_s}
\eeq
Certainly, $\M$ has dimension $\ge n-r$.
We define $$\text{vdim}(\M)=n-r$$ to be the \emph{virtual dimension} of $\M$.

The easiest case to understand is when $s$ takes values in a rank $r'$ 
subbundle $$E'\subset E$$ and is transverse to the zero section in $E'$. 
Then, $\M$ is nonsingular of dimension 
$$n-r'=\text{vdim}(M)+(r-r')\,.$$  If $E$ splits as $E=E'\oplus E/E'$, 
we can write $s=(s',0)$. We can then perturb $s$ to the section
$$
s_\epsilon=(s',\epsilon)
$$
with new zero locus given by
$$
Z(\epsilon)\subset\M\,.
$$
In particular, if $\epsilon$ can be chosen to be transverse to the zero section of $E/E'$, we obtain a smooth moduli space $Z(\epsilon)$ of the ``correct" dimension $\text{vdim}(\M)$ cut out by a transverse section $s_\epsilon$ of $E$. 
The fundamental class is
$$
[Z(\epsilon)]=c_r(E)
$$
in the (co)homology of $A$.
If we
work in the $C^\infty$ category, we can always split $E$ and pick such 
a transverse 
$C^\infty$-section.

Even when $E/E'$ has no algebraic sections (for instance if $E/E'$ is negative), the fundamental class of $Z(\epsilon)$ is clearly $c_{r-r'}(E/E')$
in the (co)homology of $\M$. 
The ``correct" moduli space, obtained
when sufficiently generic perturbations of $s$ exist or when we use $C^\infty$ sections,
has fundamental class given by the push-forward to $A$ of the top Chern class of $E/E'$.  The result is called the {\em virtual fundamental class}:
$$
[\M]^{vir}=c_{r-r'}(E/E')\in A_{\text{vdim}}(\M)\to H_{2{\text{vdim}}}(\M)\,.
$$
Here, $E/E'$, the cokernel of the derivative of the defining equations $s$,
is called the {\em obstruction bundle} of the moduli space $\M$, for reasons we explain below.

More generally $s$ need not be transverse to the zero section of any subbundle of $E$, and we must use the excess intersection
theory of Fulton-MacPherson \cite{Fu}. The limit as $t\to\infty$ of the graph of $ts$ defines a cone
$$
C_s\subset E|_\M\,.
$$
We define the virtual class to be the refined intersection 
of $C_s$ with the zero section $0_E\colon\M\into E$ inside the total space of $E$:
\beq{virdef}
[\M]^{vir}=0_E^![C_s]\in A_{\text{vdim}}(\M)\to H_{2\text{vdim}}(\M)\,.
\eeq
The result can also be expressed in terms of $c(E)s(C_s)$, 
where $c$ is the total Chern class, and $s$ is the Segre class. 

In the easy split case with $s=(s',0)$ discussed, $C_s$ is precisely $E'$. 
We recover the top Chern class $c_{r-r'}(E/E')$ of the obstruction bundle
for the virtual class.

\subsection*{\S Deformation theory}
While the basic model \eqref{model} for $\M$ rarely exists in practice
 (except in infinite dimensions), an infinitesimal version 
can be found when the moduli space admits a 
2-term deformation/obstruction theory.
The excess intersection formula \eqref{virdef} uses data only on $\M$ (rather than a neighbourhood of $\M\subset A$) and can be used in the
infinitesimal context.

At a point $p\in\M$, the basic model \eqref{model} yields
the following  exact sequence of Zariski tangent spaces
\beq{tob}
0\to T_{p\,}\M\to T_pA \Rt{ds}E_p \to\mathrm{Ob}_p\to0\,.
\eeq
So to first order, at the level of the Zariski tangent space,
 the moduli space looks like ker$\,ds$ near $p\in\M$. 
Higher order neighbourhoods of $p\in\M$ are described by 
the implicit function theorem by the zeros of the nonlinear map $\pi(s)$, where $\pi$ is the projection from $E_p$ to $\mathrm{Ob}_p$. 
The obstruction to prolonging a first order deformation of $p$ inside 
$\M$ to higher order lies in ${\mathrm{Ob}}_p$.\footnote{The 
obstruction space is not unique. Analogously, a choice of generators for the ideal of a subscheme $\M\subset A$ is not unique. For instance in our basic
 model we could have taken the obstruction bundle to be $E/E'$ \emph{or} zero. 
In each of the six approaches to curve counting,
 a natural \emph{choice} of an obstruction theory is made.}

The deformation and obstruction spaces, $T_p\M$ and ${\mathrm{Ob}}_p$, 
have dimensions differing by the virtual dimension 
$$\text{vdim}=\dim A-\rk E\ $$
and are the cohomology of a complex of \emph{vector bundles}\footnote{$B_0=TA$ and $B_1=E$ are vector bundles since $A$ is smooth and $E$ is a bundle.}
$$B_0\to B_1$$ over $\M$ restricted to $p$. 
The resolution of $T_p\M$ and ${\mathrm{Ob}}_p$ 
is the local infinitesimal method to express that $\M$ is cut out of a 
nonsingular ambient space by a section of a vector bundle. 

Li and Tian \cite{LiTianVirtual} have developed an approach
to handling deformation/obstruction theories
over $\M$.
If a global resolution $B_0\to B_1$ exists,
Li and Tian  construct a cone $C_s\subset E_1$ and intersect with the zero cycle as in \eqref{virdef} to define a virtual class on $\M$.
Due to base change issues,
the technique is difficult to state briefly, 
but the upshot is that if the deformation and obstruction spaces of a moduli problem have a difference in dimension which is \emph{constant} over $\M$ 
we can (almost always) expect a virtual cycle of the expected dimension.
\medskip

\subsection*{\S Behrend-Fantechi.}
We briefly describe a construction of the
virtual class proposed by  Behrend-Fantechi
\cite{BFNormalCone} which is equivalent and also more concise.

Dualising and globalising \eqref{tob}, we obtain the exact sequence of sheaves
$$
E^*|_\M\Rt{ds}\Omega_A|_\M\to\Omega_\M\to0\,,
$$
where the kernel of the leftmost map contains information 
about the obstructions.
The sequence factors as
$$
\spreaddiagramrows{-1pc}
\xymatrix{
E^*|_\M \rto^{ds}\dto^s & \Omega_A|_\M \ar@{=}[d] \\
I/I^2 \rto^d & \Omega_A|_\M \rto & \Omega_\M\to0\,,}
$$
where $I$ is the ideal of $\M\subset A$ and the bottom row is 
the associated exact sequence of K\"ahler differentials. 
We write $E^*|_\M\Rt{ds}\Omega_A|_\M$ as 
$$B^{-1}\to B^0,$$ 
a 2-term complex of \emph{vector bundles} because $A$ is nonsingular
 and $E$ is a bundle. The complex 
$$\{I/I^2 \to \Omega_A|_\M\}$$ is (quasi-isomorphic to) 
the  \emph{truncated cotangent complex} $\LL_\M$ of $\M$. Our data is what Behrend-Fantechi call a \emph{perfect obstruction theory}: a morphism of complexes
$$
B\udot\to \LL_\M
$$
which is an isomorphism on $h^0$ (the identity map $\Omega_\M\to\Omega_\M$) 
and a surjection on $h^{-1}$ (because $E^*\to I/I^2$ is onto). The definition can also be interpreted in terms of classical deformation theory \cite[Theorem 4.5]{BFNormalCone}.

Behrend-Fantechi show how a perfect obstruction theory leads to a cone 
in $B_1=(B^{-1})^*$ which can be intersected with the zero section to give a virtual class of dimension 
$$\text{vdim}=\rk B^0-\rk B^{-1}\,.$$
The virtual class
is the usual fundamental class when the moduli space has the correct dimension and is the top Chern class of the obstruction bundle when $\M$ is nonsingular. 
The virtual class  is also deformation invariant in an appropriate sense that would take too long to describe here.

\addtocontents{toc}{\SkipTocEntry}

\vspace{+8 pt}
\noindent Departement Mathematik \hfill Department of Mathematics \\
\noindent ETH Z\"urich \hfill  Princeton University \\
\noindent rahul@math.ethz.ch  \hfill rahulp@math.princeton.edu \\

\vspace{+8 pt}
\noindent
Department of Mathematics \\
Imperial College \\
rpwt@imperial.ac.uk

\end{document}